\documentclass[review,hidelinks,onefignum,onetabnum]{siamart220329}
\usepackage{amsfonts}
\usepackage{subcaption}
\usepackage{graphicx}
\usepackage{float}
\usepackage[dutch,english]{babel}
\usepackage{amssymb}
\usepackage{accents}
\usepackage{mathtools}
\usepackage{multicol}
\usepackage{comment}
\usepackage{caption}
\usepackage[inline]{enumitem}

\allowdisplaybreaks

\usepackage{standalone}
\usepackage{tikz}
\usetikzlibrary{shadows}
\usetikzlibrary{arrows}
\usetikzlibrary{arrows.spaced}
\usetikzlibrary{arrows.meta}
\usetikzlibrary{shapes.geometric}
\usetikzlibrary{fit}
\usetikzlibrary{backgrounds}
\usetikzlibrary{calc}

\definecolor{red}{RGB}{215,51,103}
\definecolor{blue}{RGB}{24,176,226}
\definecolor{grey}{RGB}{85,85,85}
\definecolor{mix}{RGB}{120,114,165}
\definecolor{uhhdarkgrey}{RGB}{142,169,214}
\definecolor{weirdgrey}{RGB}{105,82,64}
\definecolor{notevengrey}{RGB}{92,97,21}

\DeclareMathOperator{\syn}{Syn}

\DeclareMathOperator{\id}{Id}
\newcommand{\NN}{\mathbf{N}}

\newsiamremark{ex}{Example}
\Crefname{ex}{Example}{Examples}
\newsiamremark{remk}{Remark}
\Crefname{remk}{Remark}{Remarks}

\newcommand{\field}[1]{\mathbb{#1}}
\newcommand{\R}{\field{R}}
\newcommand{\N}{\field{N}}

\newcommand{\Z}{\field{Z}}

\DeclareMathOperator{\sgn}{sgn}

\definecolor{mycolor}{RGB}{190, 0, 190}
\definecolor{mycolor2}{RGB}{190, 80, 100}
\definecolor{wred}{rgb}{0.7,0.18,0.12}
\definecolor{wgreen}{rgb}{0.1,0.53,0.37}

\def\barray{\begin{eqnarray*}}             \def\earray{\end{eqnarray*}}
\def\beq{\begin{equation}} \def\eeq{\end{equation}}

\title{Hypernetworks: cluster synchronisation is a higher-order effect}
\date{\today}
\author{
		S\"oren von der Gracht\thanks{Department of Mathematics, Paderborn University, Germany, \href{mailto:soeren.von.der.gracht@uni-paderborn.de}{soeren.von.der.gracht@uni-paderborn.de}\funding{Partially funded by the Deutsche Forschungsgemeinschaft (DFG, German Research Foundation)--–453112019.}}
		\and
		Eddie Nijholt\thanks{Department of Mathematics, Imperial College London, United Kingdom, \href{mailto:eddie.nijholt@gmail.com}{eddie.nijholt@gmail.com}\funding{Partially supported by the Serrapilheira Institute (Grant No. Serra-1709-16124).}}
		\and
		Bob Rink\thanks{\mbox{Department of Mathematics, Vrije Universiteit Amsterdam, The Netherlands, \href{mailto:b.w.rink@vu.nl}{b.w.rink@vu.nl}}}
}
\headers{Hypernetworks: cluster synchronisation is a higher-order effect}{S\"oren von der Gracht, Eddie Nijholt,
	and Bob Rink}

\excludecomment{percent}

\newcommand{\noter}[1]{{\color{black}{}#1}}

\begin{document}
\nolinenumbers
\maketitle

\begin{abstract}
\noindent
Many networked systems are governed by non-pairwise interactions between nodes. 
The resulting higher-order interaction structure can then be encoded by means of a hypernetwork. 
In this paper we consider dynamical systems on hypernetworks by defining a class of admissible maps for every such hypernetwork. 
We explain how to classify robust cluster synchronisation patterns on hypernetworks by finding balanced partitions, and we generalise the concept of a graph fibration to the hypernetwork context. 
We also show that robust  synchronisation patterns are only fully determined by  polynomial admissible maps of high order.
This means that,  unlike in dyadic networks, cluster synchronisation on hypernetworks is a higher-order, i.e., nonlinear effect. 
We give a formula, in terms of the order of the hypernetwork, for the degree of the polynomial admissible maps that determine robust synchronisation patterns. We also demonstrate that this degree is optimal by investigating a class of examples. We conclude by displaying how this effect may cause remarkable synchrony breaking bifurcations that occur at high polynomial degree.
\end{abstract}

\section{Introduction}
\paragraph{Summary of the main results}
Recent advances in applications ranging from physics (coupled oscillator networks) over ecology (species interaction models) to social sciences (social interaction models) have indicated that, instead of by pairwise interactions, ensemble dynamics of networked real-world systems are frequently driven by simultaneous interactions of groups of network agents, so-called \emph{higher-order interactions} \cite{Ariav.2003,Levine.2017,Neuhauser.2022}. While examples of these structural aspects have been exploited in theoretical (mathematical) studies as well, a unifying framework that defines coupled dynamical systems corresponding to a higher-order hypergraph structure is still largely lacking. Consequently, this paper
\begin{itemize}
    \item generalises the concepts of coupled cell networks (\Cref{defnetwork}), admissible maps and vector fields (\Cref{def:admissible}),   graph fibrations (\Cref{def:fibration}) and quotient networks (\Cref{def:quotient}), and studies the properties of these generalisations to define, manipulate, and analyse coupled dynamical systems on hypergraphs;
    \item demonstrates that, unlike in classical (dyadic) networks,  cluster synchronisation on hypernetworks is determined by higher-order terms in the equations of motion, and is thus a purely nonlinear effect (\Cref{sec:synchrony} and \Cref{sec:examples}).
    We provide a precise expression for the polynomial degree at which cluster synchronisation is determined (\Cref{mainthr}) and show by means of examples that a lower polynomial degree is in general not sufficient (\Cref{theoremonexamplez}).
\end{itemize}
The relation between this paper and the existing concepts mentioned under the first bullet point, shall be made more precise in the background section below. Here, we would like to point out that the main result mentioned under the second bullet point distinguishes hypernetwork dynamical systems from classical (dyadic) network dynamical systems, where cluster synchronisation is known to be completely determined by the linear terms in the equations of motion \cite{Golubitsky.2005, Aguiar.2014}. As far as we know, this is one of the first examples of a dynamical phenomenon that is fundamentally different in hypernetworks than in classical networks. Furthermore, we show numerically in   \Cref{sec:examples} that this phenomenon leads to a remarkable new type of bifurcations. A systematic analytical investigation of this type of bifurcation will be presented in a separate paper. We begin by presenting an example.


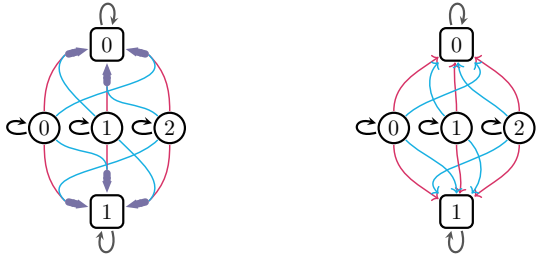
\begin{figure}[h!] 
    \centering
    \resizebox{!}{3.5cm}{
        \begin{tikzpicture}[
	square/.style = {
		regular polygon,
		regular polygon sides=4
	},
	main node/.style={
		line width=1.5pt, 
		circle,
		draw,
		font=\sffamily,
		inner sep=2pt,
		fill=white
	},
	second node/.style={
		line width=1.5pt, 
		square, 
		draw, 
		font=\sffamily,
		rounded corners, 
		inner sep=2pt
	},
	edge/.style={
		-stealth,
		shorten >=1pt,
		shorten <=1pt,
		line width=1.5pt
	},
	hyperedge/.style={
		Round Cap-{Triangle[length=3mm, width=2mm]},
		line width=5pt
	}
	]
	\node[main node, scale=1.5] (v0) at (-1.5,0) {$0$};
	\node[main node, scale=1.5] (v1) at (0,0) {$1$};
	\node[main node, scale=1.5] (v2) at (1.5,0) {$2$};
	
	\node[second node, scale=1.5] (w0) at (0,2) {$0$};
	\node[second node, scale=1.5] (w1) at (0,-2) {$1$};
	
	\path[line width=1.5pt]
	(v0) edge [edge, loop left, looseness=8] (v0)
	(v1) edge [edge, loop left, looseness=8] (v1)
	(v2) edge [edge, loop left,looseness=8] (v2)

	(w0) edge [edge, grey, loop above, looseness=8] (w0)
	(w1) edge [edge, grey, loop below, looseness=8] (w1)
	;
	
	\coordinate (he1-1) at (-1,1.75);
	\coordinate (he1-2) at (0,1);
	\coordinate (he1-3) at (1,1.75);
	\coordinate (he1-4) at (-1,-1.75);
	\coordinate (he1-6) at (0,-1);
	\coordinate (he1-5) at (1,-1.75);
	
	\path[line width=1pt]
	(v0) edge [red, in=south west, out=north] (he1-1)
	(v1) edge [blue, in=south west, out=north west] (he1-1)
	
	(v1) edge [red, in=south, out=north] (he1-2)
	(v2) edge [blue, in=south, out=north west] (he1-2)
	
	(v2) edge [red, in=south east, out=north] (he1-3)
	(v0) edge [blue, in=south east, out=north east] (he1-3)
	
	(v0) edge [red, in=north west, out=south] (he1-4)
	(v2) edge [blue, in=north west, out=south west] (he1-4)
	
	(v2) edge [red, in=north east, out=south] (he1-5)
	(v1) edge [blue, in=north east, out=south east] (he1-5)
	
	(v1) edge [red, in=north, out=south] (he1-6)
	(v0) edge [blue, in=north, out=south east] (he1-6)
	;

	\draw[hyperedge, mix] (he1-1) -- (w0);
	\draw[hyperedge, mix] (he1-2) -- (w0);
	\draw[hyperedge, mix] (he1-3) -- (w0);
	\draw[hyperedge, mix] (he1-4) -- (w1);
	\draw[hyperedge, mix] (he1-5) -- (w1);
	\draw[hyperedge, mix] (he1-6) -- (w1);
	
\end{tikzpicture}%
	}%
	\hspace{2cm}
	\resizebox{!}{3.5cm}{
        \begin{tikzpicture}[
	square/.style = {
		regular polygon,
		regular polygon sides=4
	},
	main node/.style={
		line width=1.5pt, 
		circle,
		draw,
		font=\sffamily,
		inner sep=2pt,
		fill=white
	},
	second node/.style={
		line width=1.5pt, 
		square, 
		draw, 
		font=\sffamily,
		rounded corners, 
		inner sep=2pt
	},
	edge/.style={
		-stealth,
		shorten >=1pt,
		shorten <=1pt,
		line width=1.5pt
	},
	hyperedge/.style={
		Round Cap-{Triangle[length=3mm, width=2mm]},
		line width=5pt
	}
	]
	\node[main node, scale=1.5] (v0) at (-1.5,0) {$0$};
	\node[main node, scale=1.5] (v1) at (0,0) {$1$};
	\node[main node, scale=1.5] (v2) at (1.5,0) {$2$};
	
	\node[second node, scale=1.5] (w0) at (0,2) {$0$};
	\node[second node, scale=1.5] (w1) at (0,-2) {$1$};
	
	\path[line width=1.5pt]
	(v0) edge [edge, loop left, looseness=8] (v0)
	(v1) edge [edge, loop left, looseness=8] (v1)
	(v2) edge [edge, loop left,looseness=8] (v2)

	(w0) edge [edge, grey, loop above, looseness=8] (w0)
	(w1) edge [edge, grey, loop below, looseness=8] (w1)
	;
	
	
	\path[line width=1pt]
	(v0) edge [red, in=215, out=north, ->] (w0)
	(v1) edge [blue, in=225, out=north west, ->] (w0)
	
	(v1) edge [red, in=265, out=north, ->] (w0)
	(v2) edge [blue, in=275, out=north west, ->] (w0)
	
	(v2) edge [red, in=325, out=north, ->] (w0)
	(v0) edge [blue, in=315, out=north east, ->] (w0)
	
	(v0) edge [red, in=145, out=south, ->] (w1)
	(v2) edge [blue, in=135, out=south west, ->] (w1)
	
	(v2) edge [red, in=40, out=south, ->] (w1)
	(v1) edge [blue, in=50, out=south east, ->] (w1)
	
	(v1) edge [red, in=80, out=south, ->] (w1)
	(v0) edge [blue, in=90, out=south east, ->] (w1)
	;

	
\end{tikzpicture}%
	}%
  \caption{The left figure depicts the hypernetwork that yields equations as in \eqref{firstexamplesystem}. The colouring encodes the distinct inputs of each hyperedge, while the arrows specify their targets. \noter{We define hypernetworks in Definition \ref{defnetwork}.} Replacing each hyperedge by two edges yields the  network on the right, with governing equations as in \eqref{linearbecomesclassical}. The classical network and the hypernetwork have different robust synchrony spaces.  
  }\label{fig:awesome_image3}
\end{figure}
\vspace{-4mm}
\begin{ex}
\label{ex:firstexample}
Consider the differential equations
\begin{equation}
    \label{firstexamplesystem}
\begin{split}
\dot x_0 &=  G(x_0), \,  
\dot x_1 = G(x_1), \,
\dot x_2 = G(x_2), \\
\dot y_0 &=  F(\textcolor{grey}{y_0}, \textcolor{mix}{(} \textcolor{red}{x_0}, \textcolor{blue}{x_1} \textcolor{mix}{)}, \textcolor{mix}{(} \textcolor{red}{x_1}, \textcolor{blue}{x_2} \textcolor{mix}{)}, \textcolor{mix}{(} \textcolor{red}{x_2}, \textcolor{blue}{x_0} \textcolor{mix}{)})\, , \\
\dot y_1 &=  F(\textcolor{grey}{y_1}, \textcolor{mix}{(} \textcolor{red}{x_0}, \textcolor{blue}{x_2} \textcolor{mix}{)}, \textcolor{mix}{(} \textcolor{red}{x_1}, \textcolor{blue}{x_0} \textcolor{mix}{)}, \textcolor{mix}{(} \textcolor{red}{x_2}, \textcolor{blue}{x_1} \textcolor{mix}{)})\, , \\
\end{split}
\end{equation}
for $x_0, x_1, x_2, y_0, y_1\in \mathbb{R}$. We  require that the  function $F\colon \mathbb{R} \times \mathbb{R}^2 \times \mathbb{R}^2\times \mathbb{R}^2 \to\mathbb{R}$ satisfies the invariance equations
\begin{equation}
\label{exampleF}
\begin{split}
& F(Y, (X_0, X_1), (X_2, X_3), (X_4, X_5)) =  \\
&  F(Y, (X_2, X_3), (X_0, X_1), (X_4, X_5)  ) = \\ 
&  F(Y,(X_0, X_1), (X_4, X_5), (X_2, X_3) )\, .
\end{split}
\end{equation}
In other words, the three pairs of $X$-inputs of $F$ can be exchanged without changing the value of $F$ 
 (this is why we emphasize  pairs of variables among the arguments of $F$ using additional brackets). 
We make no assumptions on the function $G\colon \R \to \R$.
 As a result of \eqref{exampleF}, we may think of the cells with states $y_0$ and $y_1$ as being targeted by three identical {\it hyperedges} of order two. The hypernetwork that encodes the  structure of equations \eqref{firstexamplesystem}     is depicted in  the left panel of \Cref{fig:awesome_image3}.

 Synchronisation occurs when groups of cells in the system evolve synchronously. Note for example that substituting $x_0=x_1=x_2$ and $y_0=y_1$ in  \eqref{firstexamplesystem} yields that $\dot x_0 = \dot x_1 = \dot x_2$ and $\dot y_0 = \dot y_1$. This implies that  the synchrony space 
$$\{x_0=x_1=x_2 \ \mbox{and} \ y_0=y_1\}$$
is invariant under the flow of any ODE of the form \eqref{firstexamplesystem}. 
The (larger) subspace
$\{y_0=y_1\}$, on the other hand,
is not. Choosing for instance 
$$F(Y, ( X_0, X_1), (X_2, X_3), (X_4, X_5) ) = X_0X_1^2 + X_2X_3^2 + X_4X_5^2$$
---which satisfies \eqref{exampleF}--- we find that 
$$
\dot y_0 =  \textcolor{red}{x_0}\textcolor{blue}{x_1}^2 + \textcolor{red}{x_1}\textcolor{blue}{x_2}^2+ \textcolor{red}{x_2}\textcolor{blue}{x_0}^2 \ \ \mbox{while}\ \ \dot y_1=   \textcolor{red}{x_0}\textcolor{blue}{x_2}^2 + \textcolor{red}{x_2}\textcolor{blue}{x_1}^2 + \textcolor{red}{x_1}\textcolor{blue}{x_0}^2 \, .
$$ 
Generically, we will thus have that
$\dot y_0 \neq  \dot y_1$ when $y_0 = y_1$, i.e., the synchrony space  $\{y_0=y_1\}$ is not flow-invariant for this $F$. 

Note that the function $F$ that we chose here is {\it nonlinear}. To see why this is important, note that
 any  {\it linear} $F$ satisfying  \eqref{exampleF} is of the form
\begin{equation}
\label{linearresponsefunction}
\begin{split}
& F(Y, (X_0, X_1), (X_2, X_3), (X_4, X_5)  )  \\  & = aY + bX_0 + cX_1 + bX_2 + cX_3 + bX_4 + cX_5\, .
\end{split}
\end{equation}
For such $F$, we see from \eqref{firstexamplesystem}   that
$$ \begin{array}{rl}
\dot y_0 = & a \textcolor{grey}{y_0} + \textcolor{mix}{(} b \textcolor{red}{x_0} + c \textcolor{blue}{x_1} \textcolor{mix}{)} + \textcolor{mix}{(} b \textcolor{red}{x_1} + c \textcolor{blue}{x_2} \textcolor{mix}{)} + \textcolor{mix}{(} b \textcolor{red}{x_2} + c \textcolor{blue}{x_0} \textcolor{mix}{)}\, , \\ 
\dot y_1 = & a \textcolor{grey}{y_1} + \textcolor{mix}{(} b \textcolor{red}{x_0} + c \textcolor{blue}{x_2} \textcolor{mix}{)} + \textcolor{mix}{(} b \textcolor{red}{x_1} + c \textcolor{blue}{x_0} \textcolor{mix}{)} + \textcolor{mix}{(} b \textcolor{red}{x_2} + c \textcolor{blue}{x_1} \textcolor{mix}{)}
 \, . 
 \end{array}
$$
It is easy to see that the right hand sides of these ODEs are equal when $y_0=y_1$, and therefore the synchrony space $\{y_0 = y_1\}$ is flow-invariant whenever $F$ is linear. 
Perhaps surprisingly, we  conclude that linear systems of the form  \eqref{firstexamplesystem} have more flow-invariant synchrony spaces than general nonlinear ones. 

One way to understand this phenomenon is to  observe from \eqref{linearresponsefunction} that any linear $F$ satisfying \eqref{exampleF} automatically satisfies the stronger invariance equations
\begin{equation}
    \label{strongerinvariance}
\begin{array}{ll}
  & F(Y, X_0, X_1, X_2, X_3, X_4, X_5)     =  \\
  & F(Y, X_2, X_1, X_0, X_3, X_4, X_5)     = \\
  & F(Y, X_0, X_1, X_4, X_3, X_2, X_5)     = \\
 &  F(Y, X_0, X_3, X_2, X_1, X_4, X_5)
 = \\
 &  F(Y, X_0, X_1, X_2, X_5, X_4, X_3)\, .
\end{array}
\end{equation}
This in turn means that any linear system of the form \eqref{firstexamplesystem} is automatically an admissible system for the (classical) coupled cell network shown in the right panel in \Cref{fig:awesome_image3}. This network  has been constructed by replacing each hyperedge of the original hypernetwork  by two edges. The admissible ODEs of this classical network are of the form
\begin{equation}
    \label{linearbecomesclassical}
\begin{split}
\dot x_0 &= G(x_0), \,
\dot x_1 = G(x_1), \,
\dot x_2 =  G(x_2), \\
\dot y_0 &= F(\textcolor{grey}{y_0}, \textcolor{red}{x_0}, \textcolor{blue}{x_1}, \textcolor{red}{x_1}, \textcolor{blue}{x_2}, \textcolor{red}{x_2}, \textcolor{blue}{x_0})\, , \\
\dot y_1 &= F(\textcolor{grey}{y_1}, \textcolor{red}{x_0}, \textcolor{blue}{x_2}, \textcolor{red}{x_1}, \textcolor{blue}{x_0}, \textcolor{red}{x_2}, \textcolor{blue}{x_1})\, , \\
\end{split}
\end{equation}
with $F$ satisfying \eqref{strongerinvariance}. One quickly checks that the synchrony space $\{y_0=y_1\}$ is invariant under the flow of all systems of the form \eqref{linearbecomesclassical}. 
\end{ex}

\paragraph{Background}
Systems of interacting dynamical units are prevalent in nature, whether it is the coordinated activity of neurons in the brain, interacting species in ecology, or opinion building in social networks. To investigate interconnected systems mathematically, one studies \emph{network dynamical systems}---coupled (nonlinear) maps or differential equations that describe individual units and their interactions. These systems behave vastly different from systems without an underlying connection structure with some of the most striking phenomena being \emph{(cluster) synchronisation}---some or all cells evolve identically---and unusual  bifurcation behaviour.

A prominent method to define dynamical systems that respect network structure is the \emph{groupoid formalism} developed by Golubitsky, Stewart, and collaborators \cite{Golubitsky.2004,Golubitsky.2006}, and Field \cite{Field.2004}. It allows to translate structural features of the network into dynamical properties using algebraic tools that can be summarized using the language of \emph{graph fibrations} \cite{DeVille.2015} and \emph{quiver representations} \cite{Nijholt.2020}. 
Key results in this field include
the classification of \emph{robust} patterns of synchrony---i.e., dynamically invariant, independent of the governing functions; compare to \Cref{ex:firstexample}---which are determined by linear systems (e.g. \cite{Aguiar.2011,Aguiar.2014,Aguiar.2018, Golubitsky.2005,Kamei.2009,Stewart.2003}),
the classification of generic bifurcations in terms of the network structure (e.g. \cite{Aguiar.2019c,Gandhi.2020,vonderGracht.2022,nijholt2019center,Soares.2018}), as well as insights into real world problems (e.g. \cite{Diekman.2013,Golubitsky.2017b,Golubitsky.1999c}) with no claim of this list being complete.



In recent years, there has been growing interest in the effect of simultaneous nonlinear interactions between three or more units---commonly referred to as \emph{higher-order interactions}---on the network dynamics. This has been ignited by developments in various disciplines: For example in neuroscience, one observes that the signal of one neuron activates or inhibits the communication channel between two other ones (cf. \cite{Ariav.2003}). In ecology, the simultaneous competition for resources of multiple species leads to nonstationary fluctuations of species abundancies typically observed in ecological networks of competing species (cf. \cite{Levine.2017}). In social sciences, multi-agent interactions can lead to a change of the average opinion in consensus dynamics (cf. \cite{Neuhauser.2022}). 
Moreover, recent results show  that higher-order interactions can emerge from data-driven model reconstruction, even when the original system is a pairwise coupled network (cf. \cite{Nijholt.2022b}). These advances suggest that also in the mathematical investigation 
the underlying structures be sharpened 
to \emph{hypernetworks} represented by \emph{hypergraphs} (\Cref{fig:awesome_image3})
. Significant progress has been made in that regard. However, most investigations have studied individual examples or specific physical systems (see for example the excellent surveys \cite{Battiston.2020,Bick.2021,Porter.2020, Torres.2021} and references therein).

From a theoretical perspective, a major obstacle to gauging the impact of higher-order interactions 
stems from the fact that the existing approaches to define network dynamical systems are either not well suited to incorporate or simply do not contain higher-order interactions at all: e.g., the groupoid formalism incorporates arbitrary group-interactions generically, while application inspired systems frequently rely on pairwise interactions only, for example by imposing additive input structure.
%
Additionally, in the existing literature it is not always clear how the higher-order interactions 
enter or shape the equations and different authors use different conventions. 
A comprehensive, unifying formalism to define admissible dynamical systems that respect the structure of a given hypernetwork is necessary. First approaches have been made only very recently and allowed for intriguing results. We want to highlight two main lines of work. One approach has been to investigate hypernetworks with an \emph{additive input structure} (cf. \cite{Aguiar.2022, Carletti.2020, Gallo.2022, Mulas.2020,Salova.2021b,Salova.2021}). On the other hand, the investigation of \emph{simplicial complexes}---i.e., hypergraphs with additional structural properties---allows for analytic results (cf. \cite{DeVille.2021b,Nijholt.2022c}). In both cases, there are tools to determine how the higher-order network structure shapes dynamics, e.g. in the form of robust cluster synchrony as well as in structure-dependent stability properties.

The goal of this paper is to take a more general stance in the sense that we consider directed hypergraphs or hypernetworks (which are more general than simplicial complexes) and define admissible maps and vector fields without the restriction to additive input structure. In particular, we generalise the groupoid formalism to hypernetworks and exploit this generalisation to characterise and understand synchronisation in the hypernetwork context. \noter{Our 
    result that cluster synchrony is a nonlinear effect further sets our construction apart from the 
    approaches to hypernetwork dynamics 
    mentioned above. While \cite{Aguiar.2022,Nijholt.2022c,Salova.2021b,Salova.2021} also prove that robust synchrony patterns are characterized by so-called balanced partitions, nonlinearity in the equations of motion is either not necessary or not further investigated in their respective formalisms.
    In fact, in \cite{Aguiar.2022} robust patterns of synchrony 
    are determined at linear degree by using the classical dyadic result, as the hypergraphs in \cite{Aguiar.2022} can be identified with a suitable bipartite dyadic graph. 
    A modification of this identification is also used in \cite{Salova.2021b,Salova.2021} to determine robust synchrony patterns in an algorithmic manner.
    References \cite{Carletti.2020, DeVille.2021b, Gallo.2022, Mulas.2020}
    mentioned above do not address cluster synchrony, but instead focus on full synchronisation.
}

\noter{We believe that the formalism presented in this paper constitutes another step towards successfully modelling real world networked systems. Our theoretical results would moreover explain, or even predict unexpected behaviour in these systems.} For example, modelling an opinion formation process according to this formalism allows to determine all robust patterns of synchrony. These, in turn, might explain the newly observed average opinion or even the occurrence of multiple opposing opinions that are shared by groups of agents.
\paragraph{Structure of the article}
This article is structured as follows. 
In \Cref{sec:hypernet}, we introduce hypernetworks and their admissible  maps as well as balanced partitions and robust synchrony subspaces. 
In \Cref{sec:fibrations}, we relate balanced partitions to hypergraph fibrations and quotient hypernetworks.
In \Cref{sec:synchrony}, we characterise robust cluster synchrony in terms of balanced partitions, and we give a polynomial degree at which cluster synchronisation is determined.  
Finally, \Cref{sec:examples} presents a class of examples that show that the polynomial degree at which cluster synchronisation is determined, found in \Cref{sec:synchrony}, is optimal. 
We conclude \Cref{sec:examples} with an example of highly unusual bifurcation behavior in a hypernetwork system, in which steady-state branches break synchrony only up to high order in the bifurcation parameter.

\paragraph{Acknowledgement}
We thank Martin Golubitsky and Ian Stewart for enlightening discussions.

\section{Hypernetworks and their admissible maps}
\label{sec:hypernet}
In this section, we formalise the idea that hypernetworks can encode the  structure of the interactions between  dynamical variables. Before introducing dynamics on hypernetworks, we first define hypernetworks as a type of directed hypergraph. Of course, the concept of the structure of a hypergraph is not new, see \cite{Ausiello.2017} for a recent survey.
\begin{definition}\label{defnetwork}
A {\it hypernetwork} is a tuple ${\rm \bf N} = (V, H, s, t)$ in which $V$ is a finite set of vertices and $H$ is a finite set of hyperedges. 
The map $s$ assigns to each hyperedge a finite ordered list of source vertices $ s(h) = (s_1(h), \ldots, s_k(h)) \in V^k$. The length $k$ of $s(h)$ is called the {\it order} of $h$, and the {\it order} of the hypernetwork is the maximum of the order of its hyperedges. The map $t: H \to V$ assigns to each hyperedge a unique target vertex. 

In addition, all vertices and hyperedges are assigned a type (chosen from some finite set), such that
\begin{enumerate}[label=\textbf{\arabic*}.]
\item if two hyperedges $h_1, h_2\in H$ have the same type, then they have equal order. Moreover, their sources $s_i(h_1)$ and $s_i(h_2)$ have the same type for each $i=1, \ldots, k$ (where $k$ is the order of $h_1$ and $h_2$), and their targets $t(h_1)$ and $t(h_2)$ have the same type; 
\item if two vertices $v_1, v_2\in V$ have the same type, then there is a type-preserving bijection $\alpha :t^{-1}(v_1)\to t^{-1}(v_2)$ between the hyperedges that target $v_1$ and $v_2$.
\end{enumerate}
A subset of vertices $V'\subset V$ such that $s(h) \in (V')^k$ for all $h\in H$ with $t(h)\in V'$ together with hyperedges $H'=\{h\in H \mid t(h)\in V'\}$ defines a \emph{sub-hypernetwork of $\NN$}, $\NN' = (V', H', s|_{H'}, t|_{H'})$. We  write $\NN' \sqsubset\NN$.
\end{definition}

\noindent Remark that a single vertex could act multiple times as a source of a hyperedge $h$, while it could also act as a source of multiple hyperedges. \noter{These properties are not standard in the literature, but they guarantee that quotients (that will be defined in Section \ref{sec:fibrations}) of hypernetworks are  hypernetworks as well.}

\Cref{defnetwork} generalises the definition of a coupled cell network as in the groupoid formalism. In fact, a coupled cell network is simply a hypernetwork in which every hyperedge has order one \noter{(cf. Def. 5.1 in \cite{Golubitsky.2006})}. Our generalisation formalises the idea that, when the hyperedges $h_1$ and $h_2$ have the same type, then the source vertices $s(h_1)=(s_1(h_1), \ldots, s_k(h_1))$ together influence the target vertex $t(h_1)$ through the hyperedge $h_1$, in exactly the same way as the source vertices $s(h_2)$ impact the target vertex $t(h_2)$ through the hyperedge $h_2$. 

\begin{ex}
  The left panel of \Cref{fig:awesome_image3} depicts a hypergraph with $5$ vertices of two types. It contains $5$ hyperedges of order $1$ of two types that each form a self-loop on one of the vertices. Additionally, it contains $6$ equal-type hyperedges of order two, with two cells of the first type as sources, and a cell of the second type as target. This example (and its generalisations) will appear at various places in this paper.
\end{ex}

\begin{figure}[ht]
    \centering
    \resizebox{.3\linewidth}{!}{
        \begin{tikzpicture}[
	square/.style = {
		regular polygon,
		regular polygon sides=4
	},
	main node/.style={
		line width=1.5pt, 
		circle,
		draw,
		font=\sffamily,
		inner sep=2pt,
		fill=white
	},
	second node/.style={
		line width=1.5pt, 
		square, 
		draw, 
		font=\sffamily,
		rounded corners, 
		inner sep=2pt
	},
	edge/.style={
		-stealth,
		shorten >=1pt,
		shorten <=1pt,
		line width=1.5pt
	},
	hyperedge/.style={
		Round Cap-{Triangle[length=3mm, width=2mm]},
		line width=5pt
	}
	]
	\node[main node, scale=1.5] (v0) at (-1.5,0) {$0$};
	\node[main node, scale=1.5] (v1) at (0,0) {$1$};
	\node[main node, scale=1.5] (v2) at (1.5,0) {$2$};
	
	\node[second node, scale=1.5] (w0) at (0,2) {$0$};
	\node[second node, scale=1.5] (w1) at (0,-2) {$1$};
	
	\path[line width=1.5pt]
	(v0) edge [edge, uhhdarkgrey, loop left, looseness=8] (v0)
	(v0) edge [edge, grey, loop left, looseness=12] (v0)
	
	(v1) edge [edge, uhhdarkgrey, out=60, in=30, looseness=6] (v1)
	(v0) edge [edge, grey] (v1)
	
	(v1) edge [edge, uhhdarkgrey] (v2)
	(v2) edge [edge, grey, loop right, looseness=8] (v2)
	;
	
	\coordinate (he1-1) at (-1,1.5);
	\coordinate (he1-2) at (0,1);
	\coordinate (he1-3) at (1,1.5);
	\coordinate (he1-4) at (-1,-1.5);
	\coordinate (he1-6) at (0,-1);
	\coordinate (he1-5) at (1,-1.5);
	
	\path[line width=1pt]
	(v0) edge [red, in=south west, out=north] (he1-1)
	(v1) edge [blue, in=south west, out=north west] (he1-1)
	
	(v1) edge [red, in=south, out=north] (he1-2)
	(v2) edge [blue, in=south, out=north west] (he1-2)
	
	(v2) edge [red, in=south east, out=north] (he1-3)
	(v0) edge [blue, in=south east, out=north east] (he1-3)
	
	(v0) edge [red, in=north west, out=south] (he1-4)
	(v2) edge [blue, in=north west, out=south west] (he1-4)
	
	(v2) edge [red, in=north east, out=south] (he1-5)
	(v1) edge [blue, in=north east, out=south east] (he1-5)
	
	(v1) edge [red, in=north, out=south] (he1-6)
	(v0) edge [blue, in=north, out=south east] (he1-6)
	;

	\draw[hyperedge, mix] (he1-1) -- (w0);
	\draw[hyperedge, mix] (he1-2) -- (w0);
	\draw[hyperedge, mix] (he1-3) -- (w0);
	\draw[hyperedge, mix] (he1-4) -- (w1);
	\draw[hyperedge, mix] (he1-5) -- (w1);
	\draw[hyperedge, mix] (he1-6) -- (w1);
	
	\begin{pgfonlayer}{background}
		\draw[rounded corners, fill=grey!20] (-2,.75) -- (2.7,.75) -- (2.7,-.75) -- (-2.7,-.75) -- (-2.7,.75) -- (-2,.75);
	\end{pgfonlayer}
\end{tikzpicture}%
	}%
    \caption{A hypernetwork that contains $3$ cells of the same type that form a classical first-order subnetwork. The $2$ additional square cells are targeted by $3$ hyperedges of order $2$ each, which only have cells of the subnetwork as sources. In contrast to \Cref{fig:awesome_image3}, we have left out the self-loop edges corresponding to the first entries of the governing functions in \eqref{runningexample}.}\label{fig:running_ex}
\end{figure}
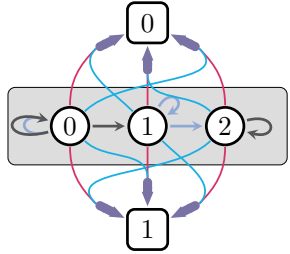

\begin{ex} \label{smallexample}
\Cref{fig:running_ex} displays a hypernetwork with $5$ vertices. Vertices $v_0, v_1$ and $v_2$ are of the same type and form a classical first-order network: they are  targeted only by edges. 
Vertices $w_0$ and $w_1$ are also of the same type and each receive one edge (from themselves; not depicted) and three equal-type hyperedges of order $2$. Both inputs of these hyperedges are taken from $v_0, v_1$ and $v_2$. 
\end{ex}
\noindent We are now ready to define hypernetwork dynamical systems in terms of the admissible maps for the hypernetwork. \noter{For a hypernetwork of order one, one recovers admissible maps for a coupled cell network (cf. Def. 6.1 in \cite{Golubitsky.2006}).}
\begin{definition}
\label{def:admissible}
Let ${\rm \bf N} = (V, H, s, t)$ be a hypernetwork.  Assume that for every $v\in V$ an {\it internal phase space} $\mathbb{R}^{n_v}$ is given in such a way that $n_{v_1}=n_{v_2}$ whenever $v_1$ and $v_2$ are of the same vertex-type. That is, vertices of the same type have identical internal phase spaces. 
A map or vector field 
$$f: \bigoplus_{v\in V} \mathbb{R}^{n_v} \to  \bigoplus_{v\in V} \mathbb{R}^{n_v}$$ 
defined on the {\it total phase space} $\bigoplus_{v\in V} \mathbb{R}^{n_v}$ 
is called ${\rm \bf N}$-{\it admissible} if it is of the form
$$f_v(x) = F_v\left(  \bigoplus_{h\, :\, t(h)=v} {\bf x}_{s(h)}  \right)\,  $$
for some {\it response function} $F_v$. Here we write 
 $${\bf x}_{s(h)} = (x_{s_1(h)}, \ldots, x_{s_k(h)}) \in \mathbb{R}^{n_{s_1(h)}} \oplus \ldots \oplus \mathbb{R}^{n_{s_k(h)}}$$ for the ordered list of source variables of $h$, and
 $$\bigoplus_{h\, :\, t(h)=v} {\bf x}_{s(h)} \in \bigoplus_{h\, :\, t(h)=v} \left( \mathbb{R}^{n_{s_1(h)}} \oplus \ldots \oplus \mathbb{R}^{n_{s_k(h)}} \right) $$
 for the list of all the input variables of $F_{v}$. 
We furthermore require that the $F_v$  satisfy the following  invariance condition: if  $\alpha: t^{-1}(v_1)\to t^{-1}(v_2)$ is any hyperedge-type-preserving bijection between the targeting hyperedges of two vertices of the same type, then 
\begin{equation}
    \label{eq:bijection_symmetry}
    F_{v_2}\left( \bigoplus_{t(h_2)=v_2} {\bf x}_{s(h_2)} \right) = F_{v_1}\left(\bigoplus_{t(h_1)=v_1} {\bf x}_{s(\alpha(h_1))} \right) \, .
\end{equation}
\end{definition}
\begin{remk}
    The invariance condition  \eqref{eq:bijection_symmetry} states that the evolution of two vertices of the same type depends ``in the same way'' on the state of the sources of any equal-type hyperedges targeting them. The condition can also be seen as a {\it local symmetry} property for the admissible vector field. 
    For example, consider the situation where $v\in V$ is targeted by two hyperedges of the same type $h_1$ and $h_2$. This implies that there is a bijection between the targeting hyperedges of $v$ that exchanges $h_1$ and $h_2$ (and keeps all other hyperedges fixed). \Cref{eq:bijection_symmetry} with $v_1 = v_2 = v$ then states that $F_v$ is invariant under the exchange of ${\bf x}_{s(h_1)}$ and ${\bf x}_{s(h_2)}$.
\end{remk}

\begin{ex} \label{smallexample_vf}
The general admissible ODE for the hypernetwork in \Cref{fig:running_ex} is
\begin{equation}\label{runningexample}
\begin{split}
    \dot x_0 & = G(x_0, \textcolor{uhhdarkgrey}{x_0}, \textcolor{grey}{x_0}),\, 
    \dot x_1 = G(x_1, \textcolor{uhhdarkgrey}{x_1}, \textcolor{grey}{x_0}), \,
    \dot x_2 = G(x_2, \textcolor{uhhdarkgrey}{x_1}, \textcolor{grey}{x_2}), \\
    \dot y_0 & = F(y_0, \textcolor{mix}{(} \textcolor{red}{x_0}, \textcolor{blue}{x_1} \textcolor{mix}{)}, \textcolor{mix}{(} \textcolor{red}{x_1}, \textcolor{blue}{x_2} \textcolor{mix}{)}, \textcolor{mix}{(} \textcolor{red}{x_2}, \textcolor{blue}{x_0} \textcolor{mix}{)})\, ,\\
    \dot y_1 & = F(y_1, \textcolor{mix}{(} \textcolor{red}{x_0}, \textcolor{blue}{x_2} \textcolor{mix}{)}, \textcolor{mix}{(} \textcolor{red}{x_1}, \textcolor{blue}{x_0} \textcolor{mix}{)}, \textcolor{mix}{(} \textcolor{red}{x_2}, \textcolor{blue}{x_1} \textcolor{mix}{)})\, .\\
    \end{split}
\end{equation}
Here, $x_0, x_1, x_2$ are the variables describing the state of the cells of the first type (circular) while $y_0, y_1$ determine the states of cells of the second type (square). The ODE is admissible precisely when 
\begin{equation}\label{exampleinvariance}
 \begin{array}{l} F(Y_0, (X_0, X_1), (X_2, X_3), (X_4, X_5)) =\\ F(Y_0, (X_2, X_3), (X_0, X_1), (X_4, X_5)) =\\ F(Y_0, (X_2, X_3), (X_4, X_5), (X_0, X_1))\, . \end{array}
 \end{equation}
This means that the three pairs of input states of $F$ can be arbitrarily permuted.
\end{ex}

\begin{figure}[ht]
    \centering
    \resizebox{\linewidth}{!}{
        \input{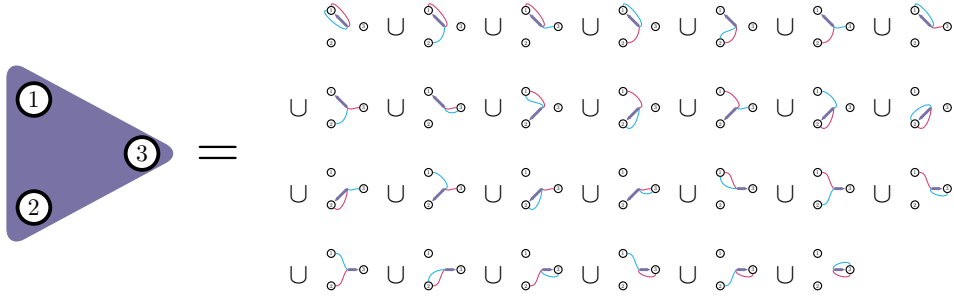}%
	}%
    \caption{An undirected hyperedge between three vertices can be regarded as the union of $27$ directed hyperedges of degree $2$. Inspired by Figure~1 in \cite{Gallo.2022}.}\label{fig:undirected}
\end{figure}

\begin{remk}
    In our definition, hyperedges are directed and target precisely one cell. However, hypergraphs that do not satisfy these conventions can often be represented as well. 
    For instance, a hyperedge with multiple targets can be considered as multiple hyperedges with one target which all have the same source vertices. 
    Similarly, an undirected hyperedge can be considered as the collection of all possible directed hyperedges connecting the involved vertices (see \Cref{fig:undirected}).
\end{remk}
    
\noindent In the remainder of this section, we describe balanced partitions and synchrony subspaces in hypernetwork systems. Consider a partition $P=\{V_1, \ldots, V_C\}$ (also called {\it colouring}) of the vertices $V$ of a hypernetwork ${\rm \bf N}$ with the property that whenever $v_1, v_2 \in V_c$ are in the same element of the partition, then they are of the same vertex-type. We can then define the polysynchrony space 
\begin{align}\label{def:SynP}
{\rm Syn}_P:= \{ x_{v_1} = x_{v_2} \ \mbox{when}\ v_1, v_2   \ \mbox{are in the same element of}\ P\, \}\, ,
\end{align}
on which the states of two vertices are synchronised when these vertices belong to the same element of the partition $P$. When the polysynchrony space is dynamically invariant for any $\NN$-admissible map, we say it is \emph{robust}.

In \Cref{def:balanced} we define what it means for a partition to be {\it balanced}. This notion can be read as the partition being compatible with the hypernetwork structure, and it is the obvious generalisation of the corresponding notion for coupled cell networks \noter{(cf. Def. 7.1 in \cite{Golubitsky.2006})}. 
\Cref{mainthr} in \Cref{sec:synchrony} below states that a partition of the cells is balanced if and only if ${\rm Syn}_P$ is robust. 

\begin{definition}\label{def:balanced}
A partition $P=\{V_1, \ldots, V_C\}$ of the vertices in a hypernetwork ${\bf N}$ is {\it balanced} if for all $v_1, v_2 \in V_c$ in the same element of the partition we have
\begin{enumerate}[label=\textbf{\arabic*}.]
\item $v_1$ and $v_2$ have the same vertex-type, i.e., $P$ is a refinement of the partition into vertex-types;
\item there is a hyperedge-type-preserving bijection $\alpha: t^{-1}(v_1) \to t^{-1}(v_2)$ such that for every hyperedge $h_1\in t^{-1}(v_1)$ of order $k$ and every source index $1\leq i\leq k$, the sources $s_i(h_1), s_i(\alpha(h_1)) \in V_d$ are also in the same element of the partition.
\end{enumerate}
\end{definition}
\begin{ex}
    \label{ex:running_synchrony}
    In \Cref{smallexample}, the partition $P=\{v_0, v_1, v_2\}\cup \{w_0, w_1\}$ of the vertices in the hypernetwork is balanced. As a matter of fact, all vertices $v_0, v_1, v_2$ receive \mbox{(hyper-)edges} from some $v_i$, which are all in the same element of the partition. Similarly, the vertices $w_0, w_1$ receive one edge from themselves, which are in the same element of the partition, and three hyperedges of order $2$ with source vertices $v_i, v_j$, which are also all in the same element of the partition. \noter{In fact, $P$ is the partition into vertex types, which is balanced for any hypernetwork.}
    
    In accordance with \Cref{mainthr} in \Cref{sec:synchrony} below, one observes from the equations  presented in \Cref{smallexample_vf} that the cluster synchrony space
    \[
     {\rm Syn}_P = \{x_0=x_1=x_2\ \mbox{and}\ y_0=y_1\}
     \]
    for this partition is flow-invariant for any admissible map, i.e., it is robust. In contrast, the (larger) space $\{y_0=y_1\}$ is not invariant under every admissible map. In fact, the partition $\{v_0\}\cup \{v_1\} \cup \{v_2\} \cup\{w_0, w_1\}$ is not balanced, as $w_0$ and $w_1$ are targeted by hyperedges whose ordered inputs come from different elements in the partition.
\end{ex}

\section{Hypergraph fibrations}
\label{sec:fibrations}
In the context of classical networks---i.e., hypernetworks of order $1$---, it is well known that the dynamics restricted to a robust synchrony subspace is that of a network as well, which is called the \emph{quotient network} (see \cite{Stewart.2003,Field.2004,Golubitsky.2005,Golubitsky.2006}). It arises by collapsing synchronous vertices to a single one and attaching arrows consistent with the original network---a construction for which it is essential that the partition is balanced. It was shown more recently, that this result is an instance of so-called \emph{graph fibrations}, which were introduced in \cite{Boldi.2002}. These are morphisms of the underlying graphs that induce linear maps sending solutions of one network dynamical system to solutions of another network dynamical system \cite{DeVille.2015}. The goal of this section is to generalise this concept to hypernetworks.

We begin by defining hypergraph fibrations (instead of hypernetwork fibrations to highlight the generalisation of graph fibrations). The definition \noter{generalises the one for classical networks, that is,  hypernetworks of order one (cf. Def. 4.1.1 in \cite{DeVille.2015} and Def. 4.2 in \cite{Nijholt.2016})}.
\begin{definition}[Hypergraph fibration]
    \label{def:fibration}
    Let $\phi\colon\NN\to\NN'$ be a map between two hypernetworks $\NN=(V,H,s,t)$ and $\NN'=(V',H',s',t')$ such that
    \begin{enumerate}[label=\textbf{\arabic*}.]
        \item $\phi$ sends vertices to vertices, i.e., $\phi(v)\in V'$ for all $v\in V$;
        \item $\phi$ sends hyperedges to hyperedges, i.e., $\phi(h)\in H'$ for all $h\in H$;
        \item $\phi$ preserves the types of vertices and hyperedges, i.e, $\phi(v)$ and $v$ as well as $\phi(h)$ and $h$ are of the same type respectively for all $v\in V$ and $h\in H$; \label{enum:type}
        \item $\phi$ sends the source vertices $s(h)\in V^k$ of a hyperedge $h\in H$ to the source vertices $s'(\phi(h))\in (V')^k$ of its image $\phi(h)\in H'$ and respects their order, i.e., $s'(\phi(h)) = (s'_1(\phi(h)), \dotsc, s'_k(\phi(h))) = (\phi(s_1(h)), \dotsc, \phi(s_k(h)))$; \label{enum:source}
        \item $\phi$ sends the unique target vertex of a hyperedge to the unique target vertex of its image, i.e., $t'(\phi(h)) = \phi(t(h))$ for all $h\in H$;
        \item and for every vertex $v\in V$, the restriction $\phi|_{t^{-1}(v)} \colon t^{-1}(v) \to (t')^{-1}(\phi(v))$ is a type-preserving bijection of hyperedges. \label{enum:fibration}
    \end{enumerate}
    Then $\phi$ is called a \emph{hypergraph fibration} or a \emph{fibration of hypernetworks}.
\end{definition}
\begin{remk}\mbox{}
    \begin{enumerate}[label=\textbf{\arabic*}.]
        \item Note that \ref{enum:source} is well-defined since hyperedges of the same type 
        have the same order. In particular, the existence of a  hypergraph fibration $\phi\colon\NN\to\NN'$ implies that the order of $\NN'$ is equal to or greater than the order of $\NN$ \noter{(note that the order of $\NN'$ may be strictly larger than that of $\NN$ if $\phi$ is not surjective.)} 
        \item Point \ref{enum:fibration} is also referred to as the \emph{fibration property} of the map $\phi$.
    \end{enumerate}
\end{remk}

\begin{figure}[h]
    \centering
    \resizebox{.3\linewidth}{!}{
        \begin{tikzpicture}[
	square/.style = {
		regular polygon,
		regular polygon sides=4
	},
	main node/.style={
		line width=1.5pt, 
		circle,
		draw,
		font=\sffamily,
		inner sep=2pt,
		fill=white
	},
	second node/.style={
		line width=1.5pt, 
		square, 
		draw, 
		font=\sffamily,
		rounded corners, 
		inner sep=2pt
	},
	edge/.style={
		-stealth,
		shorten >=1pt,
		shorten <=1pt,
		line width=1.5pt
	},
	hyperedge/.style={
		Round Cap-{Triangle[length=3mm, width=2mm]},
		line width=5pt
	}
	]
	\node[main node, scale=1.5] (v0) at (0,0) {$0$};
	
	\node[second node, scale=1.5] (w0) at (3,0) {$0$};
	
	\path[line width=1.5pt]
	(v0) edge [edge, uhhdarkgrey, loop left, looseness=8] (v0)
	(v0) edge [edge, grey, loop left, looseness=12] (v0)
	;
	
	\coordinate (he1-1) at (2,.5);
	\coordinate (he1-2) at (2,0);
	\coordinate (he1-3) at (2,-.5);
	
	\path[line width=1pt]
	(v0) edge [red, in=170, out=40] (he1-1)
	(v0) edge [blue, in=180, out=30] (he1-1)
	
	(v0) edge [red, in=170, out=10] (he1-2)
	(v0) edge [blue, in=190, out=-10] (he1-2)
	
	(v0) edge [red, in=180, out=-30] (he1-3)
	(v0) edge [blue, in=190, out=-40] (he1-3)
	;

	\draw[hyperedge, mix] (he1-1) -- (w0);
	\draw[hyperedge, mix] (he1-2) -- (w0);
	\draw[hyperedge, mix] (he1-3) -- (w0);
	
\end{tikzpicture}%
	}%
    \caption{A hypernetwork $\NN'$ that contains $2$ vertices of different types, $2$ hyperedges of order $1$ (classical edges) and $3$ hyperedges of order $2$. There is a (surjective) hypergraph fibration mapping from the hypernetwork $\NN$ in \Cref{fig:running_ex} to $\NN'$. Therefore, we call $\NN'$ a \emph{quotient} of $\NN$. As before, we have left out self-loops corresponding to self-influence of each cell.}\label{fig:running_ex_quotient}
\end{figure}
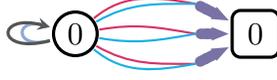
\vspace{-3mm}
\begin{ex}
    \label{ex:fibration}
    Let $\NN$ be the hypernetwork in \Cref{fig:running_ex} from our running example (\Cref{smallexample}). Furthermore, consider the second hypernetwork $\NN'$ as depicted in \Cref{fig:running_ex_quotient}. Any map $\phi\colon\NN\to\NN'$ that sends the three circular vertices $v_0,v_1$ and $v_2$ to the circular vertex $v_0$, the two square vertices $w_0$ and $w_1$ to the square vertex $w_0$, all \textsf{\textbf{\textcolor{uhhdarkgrey}{light grey}}} and \textsf{\textbf{\textcolor{grey}{grey}}} hyperedges of order $1$ to the \textsf{\textbf{\textcolor{uhhdarkgrey}{light grey}}} and \textsf{\textbf{\textcolor{grey}{grey}}} hyperedges of order $1$ respectively, the $3$ \textsf{\textbf{\textcolor{mix}{purple}}} hyperedges of order $2$ that target the square vertex $w_0$ bijectively to the $3$ \textsf{\textbf{\textcolor{mix}{purple}}} hyperedges of order $2$, and the $3$ \textsf{\textbf{\textcolor{mix}{purple}}} hyperedges of order $2$ that target the square vertex $w_1$ bijectively to the $3$ \textsf{\textbf{\textcolor{mix}{purple}}} hyperedges of order $2$, is a hypergraph fibration.
\end{ex}

\noindent The key result of \cite{DeVille.2015} \noter{(Thm. 4.3.1)} relating graph fibrations to dynamical systems translates to the hypergraph context almost immediately.
\begin{theorem}
    \label{thr:fibration}
    Let $\NN=(V,H,s,t)$ and $\NN'=(V',H',s',t')$ be two hypernetworks and $\phi\colon\NN\to\NN'$ a hypergraph fibration. Furthermore, let $f^\NN\colon \bigoplus_{v\in V} \mathbb{R}^{n_v} \to \bigoplus_{v\in V} \mathbb{R}^{n_v}$ and $f^{\NN'}\colon \bigoplus_{v'\in V'} \mathbb{R}^{n_{v'}} \to \bigoplus_{v'\in V'} \mathbb{R}^{n_{v'}}$ be admissible maps for $\NN$ and $\NN'$ respectively. That is, they are defined by
    \begin{align}
        f^\NN_v\left( \bigoplus_{w\in V} x_w \right)                &= F_v\left( \bigoplus_{h\, :\, t(h)=v} {\bf x}_{s(h)} \right) \quad \text{and} \label{eq:thr_fibration1}\\
        f^{\NN'}_{v'}\left( \bigoplus_{w'\in V'} x_{w'} \right)     &= F'_{v'}\left( \bigoplus_{h'\, :\, t'(h')=v'} {\bf x}_{s'(h')} \right), \label{eq:thr_fibration2}
    \end{align}
    where
    \[ {\bf x}_{s(h)} = \left(x_{s_1(h)}, \ldots, x_{s_k(h)}\right) \quad\text{and}\quad {\bf x}_{s'(h')} = \left(x_{s'_1(h')}, \ldots, x_{s'_{k'}(h')}\right). \]
    Assume that the internal dynamics of vertices $v\in V$ and $v'\in V'$ of the same type are governed by the same function, i.e., for every hyperedge-type-preserving bijection $\alpha \colon t^{-1}(v) \to (t')^{-1}(v')$ it holds that
    \begin{equation}
        \label{eq:fibration_condition}
        F_v\left( \bigoplus_{h\, :\, t(h)=v} {\bf x}_{s'(\alpha(h))} \right) = F'_{v'}\left( \bigoplus_{h'\, :\, t'(h')=v'} {\bf x}_{s'(h')} \right).
    \end{equation}
    Then the linear map $R_\phi \colon \bigoplus_{v'\in V'} \mathbb{R}^{n_{v'}} \to \bigoplus_{v\in V} \mathbb{R}^{n_v}$ defined by
    \begin{equation}
        \label{eq:fibration_linear}
        R_\phi \left( \bigoplus_{v'\in V'} x_{v'} \right) = \bigoplus_{v\in V} x_{\phi(v)}
    \end{equation}
    is a semiconjugacy between the admissible maps, that is, 
    \[ R_\phi \circ f^{\NN'} = f^\NN \circ R_\phi. \]
    In particular, $R_\phi$ sends solution curves of $\dot y = f^{\NN'} (y)$ on $\bigoplus_{v'\in V'} \mathbb{R}^{n_{v'}}$ to solution curves of $\dot x = f^\NN (x)$ on $\bigoplus_{v\in V} \mathbb{R}^{n_v}$.
\end{theorem}
\begin{proof}
    The proof is almost identical to that in \cite{DeVille.2015} and follows by filling in the suitable definitions and assumptions from the theorem. One computes
    \begin{align}
        \left( f^\NN \circ R_\phi \right) \left( \bigoplus_{w'\in V'} x_{w'} \right)   &= f^\NN \left( \bigoplus_{w\in V} x_{\phi(w)} \right) \label{eq:fib_proof1}\\
           &= \bigoplus_{v\in V} f^\NN_v \left( \bigoplus_{w\in V} x_{\phi(w)} \right) \label{eq:fib_proof2}\\
           &= \bigoplus_{v\in V} F_v \left( \bigoplus_{h\, :\, t(h)=v} \left( \bigoplus_{w\, : \, w\in s(h)} x_{\phi(w)} \right) \right) \label{eq:fib_proof3}\\
           &= \bigoplus_{v\in V} F_v \left( \bigoplus_{h\, :\, t(h)=v} \left( \bigoplus_{w'\, : \, w'\in s'(\phi(h))} x_{w'} \right) \right) \label{eq:fib_proof4}\\
           &= \bigoplus_{v\in V} F_v \left( \bigoplus_{h\, :\, t(h)=v} \mathbf{x}_{s'(\phi(h))} \right) \label{eq:fib_proof5}\\
           &= \bigoplus_{v\in V} F'_{\phi(v)} \left( \bigoplus_{h'\, :\, t'(h')=\phi(v)} \mathbf{x}_{s'(h')} \right) \label{eq:fib_proof6}\\
           &= R_\phi \left(\bigoplus_{v'\in V'} F'_{v'} \left( \bigoplus_{h'\, :\, t'(h')=v'} \mathbf{x}_{s'(h')} \right)\right) \label{eq:fib_proof7}\\
           &= \left(R_\phi \circ f^{\NN'}\right) \left( \bigoplus_{w'\in V'} x_{w'} \right) \label{eq:fib_proof8}
    \end{align}
    We  list the reasoning behind each equality:
    
    \begin{itemize}[noitemsep,topsep=0pt]
        \item[\eqref{eq:fib_proof1}:] by definition of $R_\phi$ \eqref{eq:fibration_linear},
        \item[\eqref{eq:fib_proof2}:] by definition of admissibility,
        \item[\eqref{eq:fib_proof3}:] by \eqref{eq:thr_fibration1},
        \item[\eqref{eq:fib_proof4}:] by \ref{enum:source} of Def.~\ref{def:fibration},
        \item[\eqref{eq:fib_proof5}:] by definition of $\mathbf{x}_{s'(h')}$,
        \item[\eqref{eq:fib_proof6}:] by \eqref{eq:fibration_condition}, since $v$ and $\phi(v)$ are of the same type and $\phi|_{t^{-1}(v)} \colon t^{-1}(v) \to (t')^{-1}(v')$ is a hyperedge-type-preserving bijection,
        \item[\eqref{eq:fib_proof7}:] by definition of $R_\phi$ \eqref{eq:fibration_linear},
        \item[\eqref{eq:fib_proof8}:] by \eqref{eq:thr_fibration2}.
    \end{itemize}
\end{proof}
\begin{remk}\mbox{}
    \begin{enumerate}[label=\textbf{\arabic*}.]
        \item Condition~\eqref{eq:fibration_condition} is well-defined. Vertices of the same type have the same internal phase space and hyperedge-type-preserving bijections $\alpha \colon t^{-1}(v) \to (t')^{-1}(v')$ for $v$ and $v'$ of the same type exist due to the local symmetry property \eqref{eq:bijection_symmetry} and the fibration property \eqref{enum:fibration} in \Cref{def:fibration} of $\phi$.
        \item We remark without further explanation that \Cref{thr:fibration} allows to extend the methods in \cite{Nijholt.2020} (to encode network structure by means of quiver representations) to hypernetworks.
    \end{enumerate}
    
\end{remk}
\begin{ex}
    We apply \Cref{thr:fibration} to the hypergraph fibration in \Cref{ex:fibration}. 
    The admissible vector field for $\NN$ is given by the right hand side of \eqref{runningexample}. Condition~\eqref{eq:fibration_condition} implies that the $\NN'$-admissible vector field is 
    $$f^{\NN'}(x_0,y_0)= \left( \begin{array}{c} G(x_0, 
    \textcolor{uhhdarkgrey}{x_0},\textcolor{grey}{x_0})  \\ F(y_0,\textcolor{mix}{(}\textcolor{red}{x_0},\textcolor{blue}{x_0}\textcolor{mix}{)},\textcolor{mix}{(}\textcolor{red}{x_0},\textcolor{blue}{x_0}\textcolor{mix}{)},\textcolor{mix}{(}\textcolor{red}{x_0},\textcolor{blue}{x_0}\textcolor{mix}{)} )\end{array}\right) \, .$$ 
    The linear map $R_\phi$ (for any choice of $\phi$) is defined as $$R_\phi(x_0,y_0)=(x_0,x_0,x_0,y_0,y_0)^T$$ 
    and it indeed semi-conjugates $f^\NN$ and $f^{\NN'}$.
\end{ex}
\noindent
Hypergraph fibrations can be used to encode structural properties of a hypernetwork (for details in the context of classical networks consult \cite{Nijholt.2020}). For example, it can readily be checked that $\NN'=(V',H',s',t')$ is (more precisely, can be identified with) a sub-hypernetwork of $\NN=(V,H,s,t)$ if and only if the inclusion $\iota\colon V'\to V$ extends to an \emph{injective} hypergraph fibration $\iota\colon\NN'\to\NN$ identifying $s'$ and $t'$ with $s$ and $t$, respectively. The linear map $R_\iota \colon \bigoplus_{v\in V}\R^{n_v} \to \bigoplus_{v\in V'} \R^{n_v}$ is then the projection
\[ R_\iota \left(\bigoplus_{v\in V}x_v \right) = \bigoplus_{v\in V'} x_{\iota(v)} = \bigoplus_{v\in V'}x_v. \]
\emph{Surjective} hypergraph fibrations are more interesting for our considerations. As in the context of classical networks, these correspond to balanced partitions \noter{(cf. Def.~5.6 in \cite{Nijholt.2020}, Lem.~5.1.1 in \cite{DeVille.2015} and Thms.~7.2 and 9.2 in \cite{Golubitsky.2006})}. 
\begin{definition}[Quotient hypernetwork]\label{def:quotient}
    Let $\mathbf{N}, \mathbf{N}'$ be two hypernetworks. If there exists a surjective hypergraph fibration $\phi\colon\mathbf{N} \to \mathbf{N}'$, then we call $\mathbf{N}'$ a \emph{quotient hypernetwork} (or simply a \emph{quotient}) of $\NN$.
\end{definition}
\begin{proposition}
    Let $\NN$ be a hypernetwork and $\NN'$ a quotient of $\NN$ corresponding to the surjective hypergraph fibration $\phi\colon\NN\to\NN'$. Then the polydiagonal
    \[ \syn_\phi := \{ x_{v_1} = x_{v_2} \text{ when } \phi(v_1)=\phi(v_2) \} \]
    is a robust synchrony subspace of $\NN$. Furthermore, any robust synchrony subspace arises in this way.
\end{proposition}
\begin{remk}
    The second part of the proposition depends on the more involved \Cref{mainthr} in \Cref{sec:synchrony} below from which we only use the result that robust synchrony implies that the corresponding partition is balanced.
\end{remk}
\begin{proof}
    The first statement follows as a corollary from \Cref{thr:fibration}: the linear map $R_\phi$ as defined in \eqref{eq:fibration_linear} is an embedding of the total phase space of $\NN'$ into the total phase space of $\NN$ whose image is the synchrony subspace $\syn_\phi$. Due to the semiconjugacy we find that $\syn_\phi$ is robust, since for any admissible map $f^{\NN}$
    \begin{multline*}
        f^{\NN}\left(\syn_\phi\right) = f^\NN\left( R_\phi\left(\bigoplus_{v'\in V'}\R^{n_{v'}}\right)\right) \\
        = R_\phi \left( f^{\NN'} \left(\bigoplus_{v'\in V'}\R^{n_{v'}}\right)\right) \subset R_\phi \left(\bigoplus_{v'\in V'}\R^{n_{v'}}\right) \subset \syn_\phi.
    \end{multline*}
    
    Conversely, assume $\syn$ is a robust synchrony subspace of $\NN$---i.e., it is defined by equality of the coordinates corresponding to certain vertices---and let $P=\{V_1,\dotsc,V_C\}$ be the corresponding partition of $V$, which is balanced due to \Cref{mainthr} below. We have that $v_1$ and $v_2$ are in the same element of the partition if and only if $x_{v_1}=x_{v_2}$ holds throughout $\syn$, and we have that $\syn=\syn_P$ \noter{-- recall that  $\rm{Syn}_P$ was defined in \eqref{def:SynP}.} We construct a quotient $\NN'$ and a surjective hypergraph fibration $\phi\colon\NN\to\NN'$. The set of vertices is given by the elements of the partition $V'=P$ and we assign the type of $v$ to $V_i$ for any $v\in V_i$---this is well-defined, since $P$ refines the partition into vertex types. Then we set $\phi(v)=V_i$ whenever $v\in V_i$ so that $\phi$ is surjective. Next, let $v_1,\dotsc,v_C$ be a set of representatives of $P$. For each $v_i$ and each $h\in t^{-1}(v_i)$ we add a hyperedge $h'$ of the same type to $H'$ such that $t'(h')=V_i$ and the inputs $s'(h')=(s'_1(h'),\dotsc,s'_k(h'))$ are defined by $s'_j(h')=V_l$, when $s_j(h)\in V_l$. 
    To define $\phi$ on hyperedges let $v\in V_i$ and $\alpha\colon t^{-1}(v) \to t^{-1}(v_i)$ be the hyperedge-type-preserving bijection from \Cref{def:balanced}. For any $h\in t^{-1}(v)$, we define $\phi(h)$ to be the hyperedge $h'$ constructed above from $\alpha(h)$, which targets the representative $v_i$. This construction is well-defined: $s_j(h)$ and $s_j(\alpha(h))$ are in the same element of $P$ for all $1\le j \le k$, since the partition is balanced. Furthermore, one may readily check that $\phi$ is indeed a hypergraph fibration.
    
    It remains to check that $\NN'$ realizes the synchrony space, i.e., that $\syn_\phi=\syn_P$. This follows immediately from the fact that, by construction, $\phi(v_1)=\phi(v_2)$ if and only if $v_1$ and $v_2$ are in the same element of the partition $P$.
\end{proof}
\noindent
In combination with \Cref{thr:fibration}, we immediately obtain
\begin{corollary}
    \label{cor:quotient}
   Let $\NN'$ be a quotient of $\NN$ under the surjective hypergraph fibration $\phi:\NN \to \NN'$. The dynamics of an $\NN$-admissible vector field restricted to the robust synchrony subspace $\syn_{\phi}$ is given by an  $\NN'$-admissible vector field. 
\end{corollary}
\begin{ex}
    Any hypergraph fibration in \Cref{ex:fibration} is surjective. Hence, $\NN'$ as depicted in \Cref{fig:running_ex_quotient} is a quotient of $\NN$ as depicted in \Cref{fig:running_ex}. It corresponds to the synchrony subspace $\{x_0=x_1=x_2 \text{ and } y_0=y_1\}$ as in \Cref{ex:running_synchrony}.
\end{ex}

\section{Balanced partitions and robust synchrony}
\label{sec:synchrony}
In the previous section we have shown that robust synchrony uniquely corresponds to surjective graph fibrations. 
In this section, we provide another characterization. 
\Cref{mainthr} is the main result regarding robust synchrony in hypernetworks, generalising the well-known result for network dynamical systems which states that balanced partitions correspond to robust synchrony. 
However, the result for hypernetworks is more subtle than the result for networks:  which synchrony spaces are robust is not determined by the linear admissible maps, but by higher order polynomial admissible maps. 
For networks (i.e. hypernetworks of order $k=1$) \Cref{mainthr} reduces to the aforementioned well-known result \noter{(cf. Thm.~7.2. in \cite{Golubitsky.2006} and Cor.~2.11 in \cite{Aguiar.2014})}. 
In what follows, we fix internal phase spaces $\R^{n_v}$ for each of the nodes.

\begin{theorem} \label{mainthr}
Let ${\bf N} = \{V, H, s, t\}$ be a hypernetwork and  $P=\{V_1, \ldots, V_C\}$  a partition of $V$ that refines the partition into vertex-types. As above, define 
$${\rm Syn}_P:= \{ x_{v_1} = x_{v_2} \ \mbox{when}\ v_1, v_2   \ \mbox{are in the same element of}\ P\ \}\, .$$ 
 The following are equivalent:
\begin{itemize}
\item[{\it i)}] The partition $P$ is balanced.
\item[{\it ii)}] ${\rm Syn}_P$ is invariant under any {\bf N}-admissible map, i.e., it is robust.
\item[{\it iii)}] ${\rm Syn}_P$ is  invariant under any polynomial {\bf N}-admissible map of degree at most $\frac{k(k+1)}{2}$, where $k$ is the order of the hypernetwork.
\end{itemize}
\end{theorem}
\noindent We split the proof of \Cref{mainthr} into two main parts. 
The first part deals with the implication ${\it i)} \implies {\it ii)}$ and is relatively straightforward and similar to the proof for dyadic networks. Note that the implication ${\it ii)} \implies {\it iii)}$ is trivial. 
The second part of the proof concerns the implication ${\it iii)} \implies {\it i)}$. This implication is considerably harder to prove than the corresponding implication for \noter{dyadic} networks.

\begin{proof}[Proof of \Cref{mainthr}, ${\it i)} \implies {\it ii)}$]
Assume that the partition $P$ is balanced and that $v_1$ and $v_2$ are in the same element of $P$. Then $v_1$ and $v_2$ are of the same vertex-type and there is a hyperedge-type-preserving  bijection $\alpha: t^{-1}(v_1) \to t^{-1}(v_2)$ so that for each hyperedge $h_1\in t^{-1}(v_1)$ and every source index $i$, the sources $s_i(h_1)$ and $s_i(\alpha(h_1))$ are also in the same element of the partition $P$. For  $x \in {\rm Syn}_P$, we thus have  that $x_{s_i(h_1)}= x_{s_i(\alpha(h_1))}$ for every $i$, which we may write as ${\bf x}_{s(h_1)} = {\bf x}_{s(\alpha(h_1))}$. It follows that for any admissible map and any $x\in {\rm Syn}_P$, we have that
\begin{align*}
f_{v_1}(x)=  & F_{v_1}\left( \bigoplus_{ t(h_1)=v_1} {\bf x}_{s(h_1)} \right)   =  F_{v_1}\left( \bigoplus_{t(h_1)=v_1} {\bf x}_{s(\alpha(h_1))} \right)   = \\ & F_{v_2}\left( \bigoplus_{t(h_2)=v_2} {\bf x}_{s(h_2)} \right) = f_{v_2}(x) \, .
\end{align*}
The third equality uses the property of an admissible map. This proves $f({\rm Syn}_P)\subset {\rm Syn}_P$, so the synchrony space is invariant under any admissible map. 
\end{proof}

\noindent Proving the implication ${\it iii)} \implies {\it i)}$ involves some combinatorics of finite number sequences and polynomials that we need to develop first. 
To this end, \noter{let $C \in \N$ be a number, and denote by $\mathcal{C} = \{1, \dots, C\}$ an ordered set of ``colours''. 
We write $\mathcal{C}^m$ for }the set of ordered sequences of length $m$ with entries in $\mathcal{C}$.
The reason for introducing sequences is because, given a partition $P = (V_1, \dots, V_C)$, we want to keep track of the elements in the partition that each hyperedge receives its inputs from. 
More precisely, to a given hyperedge $h$ of order $m$ we associate the \emph{signature} of $h$, $\mathcal{S}(h)$, as the ordered sequence $(c_1, \dots, c_m) \in \mathcal{C}^m$ where $s_{i}(h) \in V_{c_i}$ for all $i \in \{1, \dots, m\}$. We wish to understand how different monomials in a response function change when we restrict to ${\rm Syn}_P$, which is ultimately determined by the signature of each hyperedge involved.

We begin by putting a strict partial ordering $\succ$ on $\mathcal{C}^m$ as follows. Given sequences $a,b \in \mathcal{C}^m$, we set $a \succ b$ if
\begin{itemize}
\item  the number of $C$s appearing in $a$ is greater than the number of $C$s appearing in $b$, or;
\item  the number of $C$s appearing is the same for $a$ and $b$, but the number of $(C-1)$s appearing in $a$ is greater than the number of $(C-1)$s appearing in $b$, or;
\item[$\vdots$]
\item  the number of $C$s appearing is the same for $a$ and $b$, as is the number of $(C-1)$s, $(C-2)$s and so forth, up to the number of $3$s, but the number of $2$s appearing in $a$ is greater than the number of $2$s appearing in $b$.
\end{itemize}
We never have to consider the number of $1$s, as an equal number of $2$s up to $C$s means an equal number of $1$s as well (both sequences have equal length $m$). 
It follows that two sequences are only incomparable to each other if they have all symbols appearing an equal number of times. 
It is not hard to see that $\succ$ indeed defines a strict partial ordering. 
The sequence $a_C := (C, \dots, C)$ satisfies $a_C \succ b$ for all $b \not= a_C$.

Next, let $S_m$ denote the symmetric group on $m$ elements. For $a = (c_1, \dots, c_m) \in \mathcal{C}^m$ and $\sigma \in S_m$, we write $\mathcal{M}^{\sigma}_{a} \in \Z[Z_1, \dots, Z_C]$ for the monomial given by
\begin{equation*}
\mathcal{M}^{\sigma}_{a}(Z) = Z_{c_1}^{\sigma(1)}Z_{c_2}^{\sigma(2)}\dotsm Z_{c_m}^{\sigma(m)}\, .
\end{equation*}
Note that the total degree of $\mathcal{M}^{\sigma}_{a}$ is always $1+\dotsb+m = \frac{m(m+1)}{2}$.
The monomials $\mathcal{M}^{\sigma}_{a}$ will show up as the restriction to ${\rm Syn}_P$ of the terms in some conveniently chosen response functions, where $a$ will be the signature of an edge determined by $P$. 

\begin{ex}\label{ex:msigmatismsigmas}
Suppose $m=3$ and $C \geq 2$, and let $a = (2,2,1)$ and $b = (1,1,2)$. We see that $a \succ b$, as the number of $2$s appearing in $a$ is larger than the number of $2$s appearing in $b$. Let $\id \in S_3$ denote the identity permutation. We have
\begin{equation*}
\mathcal{M}^{\id}_a(Z)\! =\!  Z_2^{\id(1)}Z_2^{\id(2)}Z_1^{\id(3)}\! =\! Z_2^{1}Z_2^{2}Z_1^{3}\! =\! Z_1^3Z_2^3 \quad \text{and} \quad \mathcal{M}^{\id}_b(Z)\! =\!   Z_1^{1}Z_1^{2}Z_2^{3}\! =\! Z_1^3Z_2^3 .
\end{equation*}
Hence, we see that in this case $\mathcal{M}^{\id}_a = \mathcal{M}^{\id}_b$.
\end{ex}
\noindent Finally, we need the notion of a permutation $\tau \in S_m$ that is \emph{attuned} to a sequence $a \in \mathcal{C}^m$. 
To this end, suppose $a$ has the symbol $c$ appearing on the positions $I_c \subset \{1, \dots, m\}$, for all $c \in \mathcal{C}$.  
The permutation $\tau$ is attuned to $a$ if $\tau$ takes on its $\#I_C$ largest values on $I_C$, its next $\#I_{C-1}$ largest values on $I_{C-1}$ and so forth.

\begin{ex}
Given $a = (3,3,1,2,1) \in \mathcal{C}^5$, the permutation $\sigma \in S_5$ given by 
\begin{equation*}
\sigma(1) = 4, \sigma(2) = 5, \sigma(3) = 1, \sigma(4) = 3, \sigma(5) = 2 
\end{equation*}
is attuned to $a$. The same holds true when we switch the values of $\sigma(1)$ and $\sigma(2)$ or those of $\sigma(3)$ and $\sigma(5)$.  The identity permutation is, for instance, not attuned to $a$.
\end{ex}
\noter{
\begin{remk}
    It is not hard to see that a permutation $\tau\in S_m$ is attuned to sequence $a = (c_1, \dots, c_m)$, precisely when the rearranged sequence  $(c_{\tau(1)}, \dots, c_{\tau(m)})$ is in non-decreasing order. 
    This makes it clear that an attuned permutation always exists, and it follows that two sequences are incomparable under  $\succ$ if and only if they are the same when each is rearranged by one of its attuned permutations.
\end{remk}}
\noindent The result we need regarding these notions is the following:
\begin{lemma}\label{lemmacombbigthr}
Let $a \in \mathcal{C}^m$ be a sequence and suppose the permutation $\tau \in S_m$ is attuned to $a$. If $b \not= a$ is another sequence such that $\mathcal{M}^{\tau}_b = \mathcal{M}^{\tau}_a$, then $b \succ a$. 
\end{lemma}

\begin{proof}
Let us denote by $I^a_c, I^b_c \subset \{1, \dots, m\}$ the positions on which $a$ and $b$ have the symbol $c \in \mathcal{C}$, respectively. Note that we may write
\begin{equation}
\label{eq:degrees}
\mathcal{M}^{\tau}_a = \prod_{c=1}^C Z_c^{\sum_{i \in I^a_c} \tau(i)} \, . 
\end{equation}
By assumption, \eqref{eq:degrees} equals
\begin{equation*}
    \mathcal{M}^{\tau}_{b} = \prod_{c=1}^C Z_c^{\sum_{i \in I^b_c} \tau(i)}
\end{equation*}
so that
\begin{equation}\label{eq:relevantqeedegrees}
\sum_{i \in I^a_c} \tau(i) = \sum_{i \in I^b_c} \tau(i) \, \text{ for all } \, c \in \mathcal{C}\, .
\end{equation}
We start by looking at $c = C$. As $\tau$ is attuned to $a$, it holds that $I^a_C$ consists of the $\#I^a_C$ distinct values $i \in \{1, \dots, m\}$ for which $\tau(i)$ is largest. Hence, the only way \eqref{eq:relevantqeedegrees} can hold for $c=C$ is if either $I^a_C = I^b_C$ or $\#I^b_C > \#I^a_C$. In the latter case we indeed have $b \succ a$, whereas the former requires we look at $c = C-1$.

Suppose therefore that $I^a_c = I^b_c$ for all $c>d$, for some fixed $d \in \mathcal{C}$. Again, because $\tau$ is attuned to $a$, we see that  $I^a_d$ consists of the $\#I^a_d$ distinct values 
$$i \in \{1, \dots, m\}\setminus\left(I^a_C \sqcup \dots \sqcup I^a_{d+1}  \right) =  \{1, \dots, m\}\setminus\left(I^b_C \sqcup \dots \sqcup I^b_{d+1}  \right)$$ 
for which $\tau(i)$ is largest. \noter{Here the symbol $\sqcup$ denotes the union of disjoint sets.} The only way \eqref{eq:relevantqeedegrees} can hold for $c=d$ is when $I^a_d = I^b_d$ or $\#I^b_d > \#I^a_d$. Again, the latter case means $b \succ a$, whereas the former means we look at $c=d-1$ next. If at no point in this procedure we conclude that $b \succ a$, then eventually we arrive at $I^a_c = I^b_c$ for all $c>1$, which implies that also $I^a_1 = I^b_1$. This means $a=b$, which we excluded by assumption. Hence we indeed have $b \succ a$.
\end{proof}

\begin{proof}[Proof of \Cref{mainthr}, ${\it iii)} \implies {\it i)}$]  
To keep notation as simple as possible, we  write $h_1 \sim h_2$ to indicate that the hyperedges $h_1$ and $h_2$ are of the same hyperedge-type, and likewise use $v_1 \sim v_2$ to denote equal vertex-type for the nodes $v_1$ and $v_2$. We will write $v_1 \sim_P v_2$ to indicate that the nodes $v_1$ and $v_2$ are in the same class of the partition $P$. The set of classes or ``colours'' of $P$ will be indexed by $\mathcal{C} = \{1, \dots, C\}$.

We will prove that the partition $P$ is balanced by constructing appropriate admissible polynomial maps that allow us to count hyperedges with certain properties. We claim that $P$ is balanced if the following holds: for any hyperedge $h_0$ and sequence $a \in \mathcal{C}^m$, with $m$ the order of $h_0$, the cardinality of the set
\begin{equation}
\mathcal{N}^v_{h_0, a} := \{h \sim h_0 \mid t(h) = v \text{ and }  \mathcal{S}(h) = a\}
\end{equation}
is the same for all nodes $v$ in the same class of $P$. 
Here $\mathcal{S}(h)$ denotes the signature of $h$ determined by $P$, as defined above.
If these cardinalities match, then for any two nodes $v_1 \sim_P v_2$ we may build a bijection $\alpha: t^{-1}(v_1) \to t^{-1}(v_2)$ that maps $\mathcal{N}^{v_1}_{h_0, a}$ into $\mathcal{N}^{v_2}_{h_0, a}$ for all $a$ and $h_0$ satisfying the second condition of \Cref{def:balanced}. As $P$ refines the partition into vertex-types by assumption, this shows $P$ is indeed balanced.

We therefore fix a node $v_0$ and a hyperedge $h_0$ of order $1 \leq m \leq k$.  The response functions we will use to determine the cardinality of the sets $\mathcal{N}^v_{h_0, a}$ are as follows. Given $\sigma \in S_m$, when $v \sim v_0$ we set
\begin{equation}\label{eq:setupresponseforproof}
f_v^{\sigma}(x) = F^{\sigma}_{v}\left( \bigoplus_{ t(h)=v} {\bf x}_{s(h)} \right) = \sum_{\scriptsize{ \begin{array}{c} t(h)=v \\ h\sim h_0 \end{array} }} \hspace{-8pt} Q^{\sigma}\left({\bf x}_{s(h)}\right)\, ,
\end{equation}
in which 
\begin{equation} \label{eq:setupresponseforproof2}
Q^{\sigma}\left({\bf x}_{s(h)}\right) :=  \prod_{i = 1}^m  \left(  x_{s_i(h)}\right)^{\sigma(i)}_1 \cdot e_1^{n_v}   \, .
\end{equation}
Here $e_1^{n_v}$ denotes the first unit vector $(1,0,\dots,0)$ in $\R^{n_v}$ and $\left(  x_{s_i(h)}\right)^{\sigma(i)}_1 \in \mathbb{R}$ is the first component of $x_{s_i(h)} \in \R^{n_{s_i(h)}}$ raised to the power $\sigma(i)$.
Note that for one-dimensional internal dynamics the function \eqref{eq:setupresponseforproof2} is just given by $Q^{\sigma}\left({\bf x}_{s(h)}\right)= x_{s_1(h)}^{\sigma(1)}\dots x_{s_m(h)}^{\sigma(m)}$. 
We also point out that each $F^{\sigma}_{v}$ is polynomial of degree $1+\dotsb+m = \frac{m(m+1)}{2} \leq \frac{k(k+1)}{2}$ and satisfies the symmetry-properties imposed on response functions. 
When $v$ is not vertex-equivalent to $v_0$ we set $f^{\sigma}_v(x) = 0$. In particular, we then have that $f^{\sigma}_{v_1}$ and $f_{v_2}^{\sigma}$ agree on ${\rm Syn}_P$ whenever $v_1 \sim_P v_2$ by assumption {\it iii)}.

We now parametrize the synchrony space ${\rm Syn}_P$ by a variable $Y=(Y_c)_{c \in \mathcal{C}}$ via the map $Y\mapsto x = (x_v)_{v\in V}$ defined by $x_v:=Y_c$ if $v\in V_c$. 
It follows that for $x \in {\rm Syn}_P$ we may write ${\bf x}_{s(h)} = (Y_{c_1}, \dots, Y_{c_m})$ where $s_i(h) \in V_{c_i}$. 
Note that the indices of this latter vector are precisely the signature $\mathcal{S}(h)  = (c_1, \dots, c_m)$. 
Additionally, if we write $Z_c:= (Y_c)_1$ for all $c \in \mathcal{C}$ and $Z=(Z_c)_{c\in\mathcal{C}}$, then the first component of the vector valued function \eqref{eq:setupresponseforproof}, evaluated on ${\rm Syn}_P$, is precisely given by
\begin{equation}\label{eq:setupresponseforproof2.2}
(f^{\sigma}_v(x))_1 = \sum_{\scriptsize{ \begin{array}{c} t(h)=v \\ h\sim h_0 \end{array} }} \hspace{-8pt} \mathcal{M}^{\sigma}_{\mathcal{S}(h)}(Z)\, .
\end{equation}
Each term $h$ in the sum in \eqref{eq:setupresponseforproof2.2} yields a term $\mathcal{M}^{\sigma}_{a}$ if $\mathcal{S}(h) = a$ and we may write 
\begin{equation}\label{eq:setupresponseforproof3}
(f^{\sigma}_v(x))_1 = \sum_{a \in \mathcal{C}^m}\left(\#\mathcal{N}^v_{h_0, a}\right)\mathcal{M}^{\sigma}_{a}(Z)\, ,
\end{equation}
for $x \in {\rm Syn}_P$. It may happen that $\mathcal{M}^{\sigma}_{a} = \mathcal{M}^{\sigma}_{b}$ for distinct sequences $a,b \in \mathcal{C}^m$, as \Cref{ex:msigmatismsigmas} shows. Hence, we may not directly read off $\#\mathcal{N}^v_{h_0, a}$ as the number of monomials $\mathcal{M}^{\sigma}_{a}$ appearing in $(f^{\sigma}_v)_1|_{{\rm Syn}_P}$.

Instead we claim given a sequence $a \in \mathcal{C}^m$, if the values of $\#\mathcal{N}^v_{h_0, b}$ are known for all $b$ such that $b \succ a$, then we can  retrieve $\#\mathcal{N}^v_{h_0, a}$ from $(f^{\sigma}_v)_1|_{{\rm Syn}_P}$ for some appropriately chosen permutation $\sigma$.
To show that this indeed holds, we choose $a \in \mathcal{C}^m$ and fix a permutation $\tau \in S_m$ that is attuned to $a$. It follows from \eqref{eq:setupresponseforproof3} that the polynomial $(f^{\tau}_v)_1|_{{\rm Syn}_P}$ will involve the monomial $\mathcal{M}^{\tau}_{a}$ exactly $M_{v, a}^\tau$ times, where 
\begin{equation}\label{sumforMtau}
M_{v, a}^\tau \,\,=\hspace{-0.5cm}  \sum_{\scriptsize{ \begin{array}{c} b \in \mathcal{C}^m \\ \mathcal{M}^{\tau}_{b} = \mathcal{M}^{\tau}_{a} \end{array} }} \hspace{-0.6cm}\#\mathcal{N}^v_{h_0, b}\, .
\end{equation}
However, \Cref{lemmacombbigthr} tells us that the sum in \eqref{sumforMtau} goes only over sequences $b$ satisfying $b \succ a$, apart from $a$ itself. 
Hence, we may indeed determine $\#\mathcal{N}^v_{h_0, a}$ from the given information. 
Note that the largest sequence $a_C = (C, \dots, C)$ satisfies $\mathcal{M}^{\sigma}_{a_C}(Z) = Z_C^{m(m+1)/2}$ for any permutation $\sigma \in S_m$, whereas for all other sequences $b$ we have $\mathcal{M}^{\sigma}_{b}(Z) \not= Z_C^{m(m+1)/2}$ for all permutations $\sigma \in S_m$. 
Hence, the number $\#\mathcal{N}^v_{h_0, a_C}$ equals the coefficient in front of the monomial $\mathcal{M}^{\sigma}_{a_C}$ in $(f^{\sigma}_v)_1|_{{\rm Syn}_P}$ for any permutation $\sigma$. As any sequence  $b \not= a_C$ satisfies $a_C \succ b$, we see that we may iteratively find all numbers $\#\mathcal{N}^v_{h_0, a}$ from the maps $(f^{\sigma}_v)_1|_{{\rm Syn}_P}$ in this way.

Finally, as $f^{\sigma}_{v_1}$ and $f^{\sigma}_{v_2}$ agree on ${\rm Syn}_P$ for any $\sigma \in S_m$ whenever $v_1 \sim_P v_2$,  we see that $\#\mathcal{N}^{v_1}_{h_0, a} = \#\mathcal{N}^{v_2}_{h_0, a}$ for all such nodes $v_1, v_2$, all hyperedges $h_0$ and all signatures $a$. This shows that $P$ is indeed balanced.
\end{proof}

\noindent In \Cref{sec:examples} we present a family of examples to show that the number $\frac{k(k+1)}{2}$ is optimal and for general hypernetworks cannot be reduced. This implies in particular that {\it only} in classical networks (hypernetworks of order one) robust synchrony is determined by the linear admissible maps.

\begin{remk}
If the internal phase spaces $\R^{n_v}$ also agree for some nodes that are not of the same vertex-type, then we may define the space
$${\rm Syn}_P:= \{ x_{v_1} = x_{v_2} \ \mbox{when}\ v_1, v_2   \ \mbox{are in the same element of}\ P\ \}\, .$$
for some partition $P$ that does not refine the partition into vertex-types. 
However, such spaces are not even invariant under all constant {\bf N}-admissible maps, as can be seen by setting $f_v(x) = C_v$ for some vectors $C_v \in \R^{n_v}$ such that $C_v = C_w$ if and only the nodes $v$ and $w$ have the same vertex-type. 
\end{remk}

\begin{ex}
  Recall our running example (see \Cref{smallexample,smallexample_vf}). 
  In \Cref{ex:running_synchrony} we concluded that the partition $P=\{v_0\}\cup\{v_1\}\cup \{v_2\}\cup\{w_0, w_1\}$ is not balanced. One can check that the corresponding synchrony space $\{y_0=y_1\}$ is not robust: it is not invariant under all admissible maps, and in particular not under all polynomial admissible maps of order $\frac{2(2+1)}{2}=3$ and higher. However, one may also verify that $\{y_0=y_1\}$ is invariant under all linear and quadratic admissible maps.
\end{ex}

\section{An interesting class of examples}
\label{sec:examples}
In this section, we construct a  class of hypernetworks with interesting  ``near-synchrony'' properties. 
More precisely, given $k \geq 2$ and a hypernetwork of order at most $k$ and with $k+1$ cells of identical type,  we construct a new hypernetwork of order $k$ with a synchrony space that is not robust, but that is nevertheless invariant under every polynomial admissible map of degree strictly less than $\frac{k(k+1)}{2}$. 
These examples therefore show that the bound of $\frac{k(k+1)}{2}$ in \Cref{mainthr} cannot in general be decreased. 
We also present a brief exploration of a remarkable synchrony breaking bifurcation in one of the hypernetworks constructed in this way.  

To introduce our construction, let $S_n$ denote the symmetric group on $n$ elements (i.e., the group of permutations of $n$ elements) and write $S^{0}_n, S^{1}_n \subseteq S_n$ for the set of even and odd permutations, respectively. \noter{By $\sgn(\sigma)$ we denote the sign of a permutation $\sigma \in S_n$, which equals $0$ if $\sigma \in S_n^0$ is even, and equals $1$ if $\sigma \in S_n^1$ is odd.}

\begin{definition}\label{def:aughypernetworks}
Let $\NN$ be a given hypernetwork with $k+1 \geq 3$ nodes $v_0, \dots, v_k$, all of identical type.  
The \emph{augmented hypernetwork}  $\NN^\diamondsuit$ is obtained from $\NN$ by adding two additional nodes $w_0$ and $w_1$, their self-loops and $(k+1)!$ hyperedges of order $k$. 
The two additional nodes are of the same type, which is different from that of the $v_i$, and the additional hyperedges are likewise of a same, new type. 
These new hyperedges are labelled by the elements of the symmetric group on  $k+1$ elements, $S_{k+1}$. 
Given $\sigma \in S_{k+1}$, the hyperedge $h_{\sigma}$ satisfies 
\begin{equation}\label{augmenextrahyp}
s(h_{\sigma}) = (v_{\sigma(1)}, \dots, v_{\sigma(k)})\, , \qquad
t(h_{\sigma}) = w_{\sgn(\sigma)}\, , 
\end{equation}
where $S_{k+1}$ is understood to act on the ordered set $\{0, \dots, k\}$ and $v_{\sigma(0)}$ is therefore the only node in $\NN$ that is not a source of $h_{\sigma}$.
We will refer to $\NN$ as the \emph{core} of the augmented hypernetwork $\NN^\diamondsuit$. 
\end{definition}
\noindent The order of the augmented hypernetwork in \Cref{def:aughypernetworks} is $k$, provided the core has order $k$ or less. 
In particular, this holds when the core is a classical (dyadic) network.
\begin{ex}\label{ex:augdisconncore}
    Let $\NN$ be the classical network consisting of three disconnected nodes, with only a single self-loop for each of them. 
    It follows that the augmented hypernetwork $\NN^\diamondsuit$ is the one shown in the left panel of \Cref{fig:awesome_image3}. 
    Here the circular nodes belong to $\NN$, whereas the square ones are the newly added $w_0$ and $w_1$.
    Note that the core $\NN$ is not required to be connected in \Cref{def:aughypernetworks}.
\end{ex}

\begin{ex}\label{ex:augrunnexam}
The hypernetwork from \Cref{smallexample} is the augmented hypernetwork with core the classical three-node network shown within the box in \Cref{fig:running_ex} (i.e. consisting of the circular nodes and the arrows in between).
\end{ex}
\noindent Given an augmented network $\NN^\diamondsuit$ with core $\NN$, we shall  denote by $x_0, \ldots, x_k$ the dynamical variables of the cells $v_0, \ldots, v_k$ and by $y_0, y_1$ the variables of the cells $w_0, w_1$. 
For convenience, we assume from here on out that all cells have a one-dimensional internal phase space. 
It follows that the equations of motion for the $y$-variables are given by
\begin{align}\label{Dsystems}
\begin{array}{l} 
\dot y_0 =  F\left(y_0,   \bigoplus_{\sigma \in S^0_{k+1}} {\bf x}_{\sigma} \right)\, , \, \,
\dot y_1 =  F\left(y_1,   \bigoplus_{\sigma \in S^1_{k+1}} {\bf x}_{\sigma} \right)\, .
\end{array}
\end{align}
Here we used the notation
\begin{equation*}
    {\bf x}_{\sigma} = (x_{\sigma(1)}, \dots x_{\sigma(k)}) \in \mathbb{R}^k\, 
\end{equation*}
for the source variables of the hyperedge $h_{\sigma}$. The response function
\begin{equation*}
    F = F\left(Y,  \bigoplus_{\sigma \in S^0_{k+1}}  {\bf X}_{\sigma}\right) \ \mbox{from}\  \mathbb{R} \oplus \!\! \bigoplus_{\sigma \in S^0_{k+1}} \mathbb{R}^k \ \mbox{to}\ \mathbb{R}
\end{equation*}
is assumed to be invariant under any permutation of the $\frac{(k+1)!}{2}$ entries ${\bf X}_{\sigma}\in \mathbb{R}^k$, which implies in particular that the notation  
$ \bigoplus_{\sigma \in S^{0}_{k+1}} {\bf X}_{\sigma}$ for the arguments of $F$ is unambiguous: the ${\bf X}_{\sigma}$ can be substituted into $F$ in arbitrary order. 
As $S^0_{k+1}$ and $S^1_{k+1}$ have the same cardinality, the invariance of $F$ likewise means that the expression for $\dot{y}_1$ in \Cref{Dsystems} is well-defined. 
Note that each $\dot{x}_i$ depends only on the $x$-variables, according to the hypernetwork structure of the core $\NN$. This is because the core is a sub-hypernetwork of the augmented hypernetwork. 

Our main result about these augmented hypernetworks is the following.
\begin{theorem}\label{theoremonexamplez}
Let $\NN^{\diamondsuit}$ be an augmented hypernetwork whose core consists of $k+1$ nodes. 
Assume one-dimensional internal dynamics for each of the nodes, and write $F$ for the response function of the $y$-nodes as in \Cref{Dsystems}.
The space $\{y_0 = y_1\}$ is invariant for all admissible systems for $\NN^{\diamondsuit}$ with $F$ polynomial of total degree strictly less than $\frac{k(k+1)}{2}$, but not for all polynomials $F$ of total degree  $\frac{k(k+1)}{2}$. 
In particular, $\{y_0 = y_1\}$ is not a robust synchrony space.
\end{theorem}
\noindent
As the dynamics of the cells in the core does not depend on the $y$-variables, we see that invariance of the space $\{y_0 = y_1\}$ is equivalent to the condition 
\begin{align}\label{eq:y0isy1equi}
F\left(y,   \bigoplus_{\sigma \in S^0_{k+1}} {\bf x}_{\sigma} \right) =  F\left(y,   \bigoplus_{\sigma \in S^1_{k+1}} {\bf x}_{\sigma} \right)
\end{align}
for all $y = y_0 = y_1 \in \R$ and $x = (x_0, \dots, x_k) \in \R^{k+1}$. 
This explains why the result of \Cref{theoremonexamplez} does not depend on the core $\NN$.
At the end of this section we will investigate a phenomenon in augmented hypernetworks that does depend on specifics of the core. 
\Cref{eq:y0isy1equi} also suggests that in order to prove \Cref{theoremonexamplez}, we first need to gather results on functions with the symmetry properties of $F$. To this end we have the following lemmas.

 \begin{lemma}\label{lem:formain0}
 Let $Q\colon \bigoplus_{\sigma \in S_{k+1}^0} \R^k \rightarrow \R$ be a function that is invariant under all permutations of its $\#S_{k+1}^0$ entries in $\R^k$. 
 Then we have
   \begin{equation*}
        Q\left(\bigoplus_{\sigma \in S^0_{k+1}} {\bf x}_{\sigma} \right) = Q\left(\bigoplus_{\sigma \in S^1_{k+1}} {\bf x}_{\sigma} \right) 
    \end{equation*}
 for $x = (x_0, \dots, x_k) \in \R^{k+1}$ whenever $x_i = x_j$ for some distinct $i,j \in \{0, \dots, k\}$.
\end{lemma}

\begin{proof}
    Let $i$ and $j$ be as in the lemma and let $\kappa \in S_{k+1}^1$ denote the transposition that interchanges $i$ and $j$, while leaving the other indices fixed. 
    By assumption we have $x_{\kappa(\ell)} = x_{\ell}$ for all $\ell \in \{0, \dots, k\}$. 
    It follows that for all $\sigma \in S_{k+1}$ and $m \in \{1, \dots, k\}$ it holds that
    \begin{equation*}
        ({\bf x}_{\kappa\sigma})_{m} = x_{\kappa(\sigma(m))} = x_{\kappa(\ell)} = x_{\ell} = x_{\sigma(m)} = ({\bf x}_{\sigma})_{m}\, ,
    \end{equation*}
    where we have set $\ell = \sigma(m)$. 
    Hence, we see that ${\bf x}_{\kappa\sigma} = {\bf x}_{\sigma}$ for all $\sigma \in S_{k+1}$. 
    From the fact that $\kappa S^0_{k+1} = S^1_{k+1}$ as sets, together with the symmetry properties of $Q$, we indeed find
    \begin{equation*}
        Q\left(\bigoplus_{\sigma \in S^0_{k+1}} {\bf x}_{\sigma} \right) = Q\left(\bigoplus_{\sigma \in S^0_{k+1}} {\bf x}_{\kappa\sigma} \right) = Q\left(\bigoplus_{\sigma \in S^1_{k+1}} {\bf x}_{\sigma} \right) \, ,
    \end{equation*}
    which completes the proof.
\end{proof}

\begin{lemma}\label{lem:formain1}
    Let $Q$ be as in \Cref{lem:formain0} and assume in addition that this function is polynomial.   
    There exists a polynomial $S \colon \R^{k+1} \rightarrow \R$ such that
  \begin{equation}\label{eq:factorisation}
        Q\left(\bigoplus_{\sigma \in S^0_{k+1}} {\bf x}_{\sigma} \right) - Q\left(\bigoplus_{\sigma \in S^1_{k+1}} {\bf x}_{\sigma} \right) = S(x) \prod_{\substack{i,j=0 \\ i > j }}^k (x_i - x_j)
    \end{equation}
    for all $x = (x_0, \dots, x_k) \in \R^{k+1}$.
 \end{lemma}

 \begin{proof}
   It follows from \Cref{lem:formain0} that the left hand side of \eqref{eq:factorisation} vanishes whenever $x_i = x_j$ for some distinct $i,j \in \{0, \dots, k\}$. 
   Let us fix two such indices $i \not= j$.
   Any polynomial $P$ in the variables $x = (x_0, \dots, x_k)$ may be written as 
\begin{equation*}
    P(x) = (x_i - x_j)T(x) + R(x_0, \dots, \hat{x}_i, \dots, x_k),  
\end{equation*}
    for some polynomials $T$ and $R$, and where $\hat{x}_i$ means that $R$ does not depend on $x_i$. 
    This can be seen by setting $x_i = (x_i - x_j) + x_j$ and separating out multiples of $x_i - x_j$. 
    If $P$ vanishes when we set $x_i = x_j$ then necessarily $R = 0$, so that $x_i - x_j$ divides $P$.
    
    Returning to \eqref{eq:factorisation}, we conclude that the left hand side is divisible by $x_i - x_j$ for all pairs of distinct indices $(i,j)$. 
    If we impose $i>j$ then the factors $x_i - x_j$ are all different irreducible polynomials (i.e., not differing by a unit). 
    Using that the polynomial ring $\R[x_0, \dots, x_k]$ is a unique factorization domain (a UFD), we conclude that the left hand side of \eqref{eq:factorisation} is indeed divisible by 
    \begin{equation*}
    \prod_{\substack{i,j=0 \\ i > j }}^k (x_i - x_j)\, ,
   \end{equation*}
    from which \eqref{eq:factorisation} follows.
  \end{proof}

\begin{corollary}\label{cor:formain0}
    Let $Q$ be a polynomial as in \Cref{lem:formain0,lem:formain1}. If $Q$ is of total degree less than $\frac{k(k+1)}{2}$, then 
     \begin{equation}\label{eq:equalforlowdegreee}
        Q\left(\bigoplus_{\sigma \in S^0_{k+1}} {\bf x}_{\sigma} \right) = Q\left(\bigoplus_{\sigma \in S^1_{k+1}} {\bf x}_{\sigma} \right) 
    \end{equation}
    for all $x = (x_0, \dots, x_k) \in \R^{k+1}$.
\end{corollary}

\begin{proof}
From \Cref{lem:formain1} we get \eqref{eq:factorisation} for some polynomial $S$. The left hand side of \eqref{eq:factorisation} has degree less than $\frac{k(k+1)}{2}$, whereas the factor $\prod_{\substack{i,j=0 \\ i > j }}^k (x_i - x_j)$ on the right hand side has degree  $\frac{k(k+1)}{2}$. So \eqref{eq:factorisation} can only hold when $S = 0$, which proves  \eqref{eq:equalforlowdegreee}.
\end{proof}
\noindent Finally, we introduce a symmetric function for the proof of \Cref{theoremonexamplez}. Given $k \in \N$, we define the \emph{power sum symmetric polynomial} $P_{(k)} \colon \bigoplus_{\sigma \in S^0_{k+1}} \mathbb{R}^k \to \mathbb{R}$ as
\begin{equation}\label{eq:powersumsym1}
P_{(k)}\left(\bigoplus_{\sigma \in S^0_{k+1}}  {\bf X}_{\sigma}\right) = \sum_{\sigma \in S^0_{k+1}}X_{\sigma,1}^{1}X_{\sigma,2}^{2}\dotsm X_{\sigma,k}^{k}\, ,
\end{equation}
where ${\bf{X}}_{\sigma} = (X_{\sigma,1}, \dots, X_{\sigma,k}) \in \R^k$ for $\sigma \in S^0_{k+1}$.
Note that $P_{(k)}$ is symmetric under all permutations of its $\#S_{k+1}^0$ entries ${\bf{X}}_{\sigma}$ in $\R^k$ and has degree $1+ \dotsb+k = \frac{k(k+1)}{2}$.

\begin{proof}[Proof of \Cref{theoremonexamplez}]
We first show that the space $\{y_0 = y_1\}$ is not in general invariant when the response function $F$ in \Cref{Dsystems} is a polynomial of degree $\frac{k(k+1)}{2}$. 
To this end, we set 
\begin{equation*}
F\left(Y,   \bigoplus_{\sigma \in S^0_{k+1}} {\bf X}_{\sigma} \right) = P_{(k)}\left(\bigoplus_{\sigma \in S^0_{k+1}} {\bf X}_{\sigma} \right)\, ,
\end{equation*}
where $P_{(k)}$ is defined by \eqref{eq:powersumsym1}. It follows that
\begin{equation}\label{fwithspecialPeven}
   \dot y_0 =  F\left(y_0,   \bigoplus_{\sigma \in S^0_{k+1}} {\bf x}_{\sigma} \right) = \!\! \sum_{\sigma \in S^0_{k+1}} x_{\sigma(1)}^1\dotsm x_{\sigma(k)}^k \, ,
\end{equation}
and similarly
\begin{equation}\label{fwithspecialPodd}
  \dot y_1 =   F\left(y_1,   \bigoplus_{\sigma \in S^1_{k+1}} {\bf x}_{\sigma} \right) = \!\! \sum_{\sigma \in S^1_{k+1}} x_{\sigma(1)}^1\dotsm x_{\sigma(k)}^k\, .
\end{equation}
It is not hard to see that the right hand sides of \Cref{fwithspecialPeven,fwithspecialPodd} are not equal. For example, we can only have
\[x_{\sigma(1)}^1\dotsm x_{\sigma(k)}^k = x_0^0x_1^1\dotsm x_k^k \]
if $\sigma \in S_{k+1}$ is the identity. Hence, this particular monomial appears in  \Cref{fwithspecialPeven} but not in  
 \Cref{fwithspecialPodd}. This shows that  the polynomials on the right hand sides of both equations are indeed different, and hence that $\{y_0=y_1\}$ is not robust.

On the other hand,  suppose now that $F$ is a polynomial of degree $d < \frac{k(k+1)}{2}$ satisfying the required symmetry conditions. It follows that we may write
\begin{equation*}
F\left(Y,   \bigoplus_{\sigma \in S^0_{k+1}} {\bf X}_{\sigma} \right) = \sum_{\ell = 0}^{d}\,  Y^{\ell}\, Q_{\ell}\left(\bigoplus_{\sigma \in S^0_{k+1}} {\bf X}_{\sigma} \right)\, ,
\end{equation*}
where each polynomial $Q_{\ell}$ is invariant under all permutations of its $\frac{(k+1)!}{2}$ entries ${\bf X}_{\sigma}$, and of degree strictly less than $\frac{k(k+1)}{2}$. 
By \Cref{cor:formain0} we have
\begin{equation}\label{equalitytoshowsym}
    Q_{\ell}\left(\bigoplus_{\sigma \in S^0_{k+1}} {\bf x}_{\sigma} \right) = Q_{\ell}\left(\bigoplus_{\sigma \in S^1_{k+1}} {\bf x}_{\sigma} \right)
\end{equation}
for all $\ell \in \{0, \dots, d\}$. This in turn implies that
\begin{equation*}
\dot{y}_0 = F\left(y_0,   \bigoplus_{\sigma \in S^0_{k+1}} {\bf x}_{\sigma} \right) = F\left(y_1,   \bigoplus_{\sigma \in S^1_{k+1}} {\bf x}_{\sigma} \right) = \dot{y}_1
\end{equation*}
whenever $y_0 = y_1$, which proves that the space $\{y_0 = y_1\}$ is dynamically invariant.
\end{proof}

\begin{remk}\label{symmetryrem}
The synchrony space $\{y_0 = y_1\}$ is the fixed point space of the map
$$S: (x_0, \ldots, x_k, y_0, y_1) \mapsto (x_{0}, \ldots, x_{k}, y_{1}, y_0 )\ .$$ 
From the proof of \Cref{theoremonexamplez} it is clear that $S$ is a symmetry of the admissible vector field for the augmented hypernetwork whenever $F$ is polynomial of degree less than $\frac{k(k+1)}{2}$, but not in general. 
As fixed point spaces are dynamically invariant for equivariant systems \noter{(see for example Theorem~8.4 in \cite{Golubitsky.2023})}, this offers an alternative interpretation of \Cref{theoremonexamplez}.
\end{remk}
\noindent We conclude this section by describing a type of bifurcation that appears to occur abundantly in  augmented hypernetworks. 
The idea is that, even though $\{y_0 = y_1\}$ is not a robust synchrony space, its invariance under polynomial response functions of sufficiently low degree, makes it act as a ``ghost'' synchrony space that can  still influence bifurcations.
In particular, in various examples we found steady-state bifurcation branches $(x(\lambda), y(\lambda))$ in which $y_0(\lambda)$ and $y_1(\lambda)$  are not exactly equal, but do agree up to unusually high degree in the bifurcation parameter $\lambda$. 
We have dubbed this phenomenon ``reluctant synchrony breaking'' and we investigate it in detail in a companion paper. 
Here we only illustrate it  in one numerical example.

\begin{ex}\label{ex:leadingreluctant00}
We revisit the augmented hypernetwork from \Cref{smallexample,ex:augrunnexam}, of which the admissible ODEs are given by \Cref{runningexample} in \Cref{smallexample_vf}. 
Instead of a single admissible vector field, we consider a family of them by defining for each $\lambda \in \R$ the response functions
\begin{align}
   \label{num:res-g} &G_{\lambda}(X_0,X_1,X_2) =  -X_0+X_1-X_2+8\lambda X_0+4X_0^2 \, \text{ and }\\
    \label{num:res-f} &F_{\lambda}(Y, (X_0, X_1), (X_2, X_3), (X_4, X_5) ) =  -5Y + 14 \lambda\\ \nonumber
    &-h(10X_0-12X_1)-h(10X_2-12X_3)-h(10X_4-12X_5)\, ,   
\end{align}
in which
\begin{equation}\label{num:res-h} 
   h(x) = \sin(x)+\cos(x)-1 = \sqrt{2}\sin(x+\frac{\pi}{4}) - 1\, .
\end{equation}
For every  value of $\lambda$, the function $F_{\lambda}$ satisfies the  symmetry conditions required in  \Cref{exampleinvariance}. 
As $G_{0}(0,0,0) = F_{0}(0, (0, 0), (0, 0), (0, 0) ) = 0$, the resulting admissible vector field  has a fixed point at the origin for $\lambda=0$. 
 Moreover, by construction the Jacobian of the admissible vector field around the origin is singular, so that we may expect a steady-state bifurcation to occur as $\lambda$ is varied near $0$.

\Cref{firstnum} shows the results of a numerical bifurcation analysis of the problem, in which we found a stable branch of steady states for small $\lambda<0$ and a stable branch of steady states for small $\lambda >0$. 
\Cref{firstnum-a} is a bifurcation diagram showing the values of the different components $x_i(\lambda)$ and $ y_i(\lambda)$ on these stable branches. 
The branch for $\lambda<0$  appears to lie in the robust synchrony space $\{x_0=x_1=x_2 \ \mbox{and} \ y_0=y_1\}$. \noter{At $\lambda = 0$, the core (corresponding to the $x$-values) undergoes a synchrony breaking bifurcation, causing the $x_i$ to become fully non-synchronous for $\lambda>0$. 
For small positive values of $\lambda$, the $y$-values on this branch seem to remain equal, but at higher values of $\lambda$, it becomes clear that a small separation occurs.  
\Cref{firstnum-b} corroborates this observation, showing a non-linear departure from the space $\{y_0 = y_1\}$, with \Cref{firstnum-c} indicating that in fact $y_0(\lambda) - y_1(\lambda) \sim \lambda^3$ for $\lambda >0$.}

\begin{figure}[h!]
\centering
\begin{subfigure}{.65\textwidth}
\centering
\includegraphics[width=\linewidth]{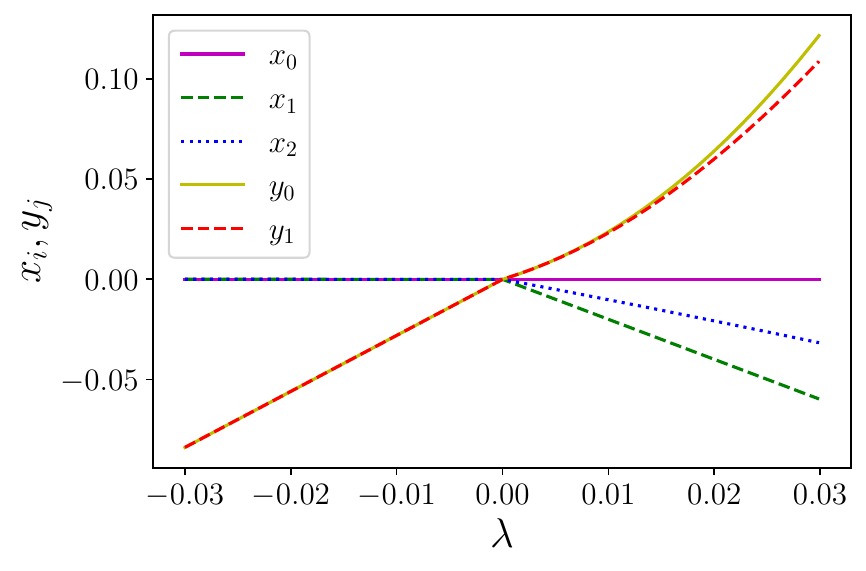}
\caption{The stable branches of a synchrony breaking bifurcation.}
\label{firstnum-a}
\end{subfigure}\hfill\par\medskip
\begin{subfigure}{0.45\textwidth}
\centering
\includegraphics[width=\linewidth]{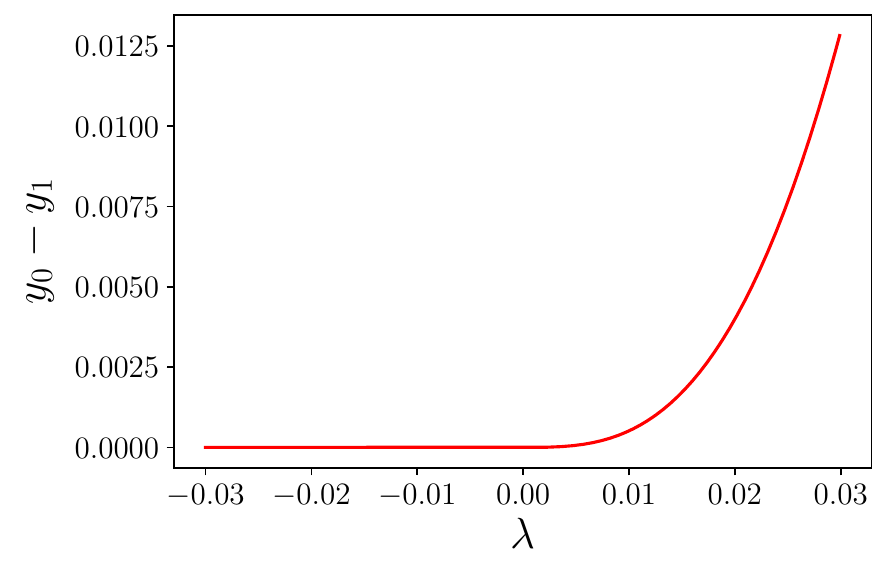}
\caption{The difference between the $y$-nodes along the stable branches.}
\label{firstnum-b}
\end{subfigure}\hfill 
\begin{subfigure}{0.45\textwidth}
\centering
\includegraphics[width=\linewidth]{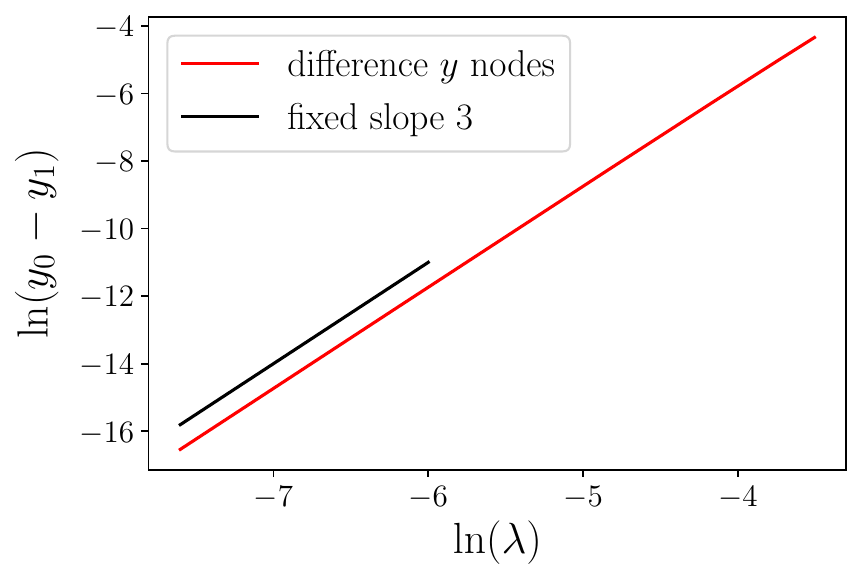}
\caption{A log–log plot of the difference between the $y$-nodes, for $\lambda > 0$.}
\label{firstnum-c}
\end{subfigure}%
\caption{Numerically obtained bifurcation diagram for a family of systems corresponding to the augmented hypernetwork shown in \Cref{fig:running_ex}. For comparison, the black line segment in the log-log plot has slope $3$, indicating that $y_0(\lambda) - y_1(\lambda) \sim \lambda^3$. }
\label{firstnum}
\end{figure}
 \noter{Such ``reluctant synchrony breaking'' would be highly anomalous in bifurcations of general vector fields. However, we  claim (and do not prove here but elsewhere) that it happens in generic one-parameter families of hypernetwork vector fields of the form \eqref{runningexample}. In particular, it is not an artifact of our particular choice of response functions.} Rather, our response functions \noter{\eqref{num:res-g} and \eqref{num:res-f}} are merely chosen to guarantee stability of each branch for the correct sign of $\lambda$, to produce clear pictures, and to illustrate that the reluctant behavior $y_0(\lambda) - y_1(\lambda) \sim \lambda^3$ is not just a consequence of (\noter{low-degree}) polynomial response functions.

\Cref{firstnum} was obtained by forward integrating the equations of motion for each of $600$ equidistributed values of $\lambda \in [-0.03, 0.03]$ \noter{(in $[0.0005, 0.03]$ for \Cref{firstnum-c})}, using Euler's method with time steps of $0.1$. 
For each value of $\lambda$ integration was performed up to $t=2000$ and starting from the point $(x_0, x_1, x_2, y_0, y_1) = (0.1, -0.2, 0.3, 0.4, 0.5)$ in phase space. Note that only stable branches can be visualised in this way.
\end{ex}

\bibliographystyle{siamplain}
\bibliography{bibliography}

\providecommand{\noopsort}[1]{}
\begin{thebibliography}{10}

\bibitem{Aguiar.2022}
{\sc M.~Aguiar, C.~Bick, and A.~Dias}, {\em Network dynamics with higher-order
  interactions: Coupled cell hypernetworks for identical cells and synchrony},
  Nonlinearity, 36 (2023), p.~4641,
  \url{https://doi.org/10.1088/1361-6544/ace39f}.

\bibitem{Aguiar.2019c}
{\sc M.~A. Aguiar, A.~P. Dias, and P.~Soares}, {\em The steady-state lifting
  bifurcation problem associated with the valency on networks}, Physica D:
  Nonlinear Phenomena, 390 (2019), pp.~36--46, \url{https://doi.org/10.1016/j.
  physd.2018.10.006}.

\bibitem{Aguiar.2011}
{\sc M.~A.~D. Aguiar, P.~Ashwin, A.~P.~S. Dias, and M.~J. Field}, {\em Dynamics
  of coupled cell networks: Synchrony, heteroclinic cycles and inflation},
  Journal of Nonlinear Science, 21 (2011), pp.~271--323,
  \url{https://doi.org/10.1007/s00332-010-9083-9}.

\bibitem{Aguiar.2014}
{\sc M.~A.~D. Aguiar and A.~P.~S. Dias}, {\em The lattice of synchrony
  subspaces of a coupled cell network: Characterization and computation
  algorithm}, Journal of Nonlinear Science, 24 (2014), pp.~949--996,
  \url{https://doi.org/10.1007/s00332-014-9209-6}.

\bibitem{Aguiar.2018}
{\sc M.~A.~D. Aguiar and A.~P.~S. Dias}, {\em An overview of synchrony in
  coupled cell networks}, in Modeling, Dynamics, Optimization and Bioeconomics
  III, A.~A. Pinto and D.~Zilberman, eds., vol.~224 of Springer Proceedings in
  Mathematics {\&} Statistics, {Springer International Publishing}, Cham, 2018,
  pp.~25--48, \url{https://doi.org/10.1007/978-3-319-74086-7_2}.

\bibitem{Ariav.2003}
{\sc G.~Ariav, A.~Polsky, and J.~Schiller}, {\em Submillisecond precision of
  the input-output transformation function mediated by fast sodium dendritic
  spikes in basal dendrites of ca1 pyramidal neurons}, Journal of Neuroscience,
  23 (2003), pp.~7750--7758,
  \url{https://doi.org/10.1523/JNEUROSCI.23-21-07750.2003}.

\bibitem{Ausiello.2017}
{\sc G.~Ausiello and L.~Laura}, {\em Directed hypergraphs: Introduction and
  fundamental algorithms---a survey}, Theoretical Computer Science, 658 (2017),
  pp.~293--306, \url{https://doi.org/10.1016/j. tcs.2016.03.016}.

\bibitem{Battiston.2020}
{\sc F.~Battiston, G.~Cencetti, I.~Iacopini, V.~Latora, M.~Lucas, A.~Patania,
  J.-G. Young, and G.~Petri}, {\em Networks beyond pairwise interactions:
  Structure and dynamics}, Physics Reports,  (2020), pp.~1--92,
  \url{https://doi.org/10.1016/j. physrep.2020.05.004}.

\bibitem{Bick.2021}
{\sc C.~Bick, E.~Gross, H.~A. Harrington, and M.~T. Schaub}, {\em What are
  higher-order networks?}, July 2022, \url{https://arxiv.org/abs/2104.11329}.

\bibitem{Boldi.2002}
{\sc P.~Boldi and S.~Vigna}, {\em Fibrations of graphs}, Discrete Mathematics,
  243 (2002), pp.~21--66, \url{https://doi.org/10.1016/S0012-365X(00)00455-6}.

\bibitem{Carletti.2020}
{\sc T.~Carletti, D.~Fanelli, and S.~Nicoletti}, {\em Dynamical systems on
  hypergraphs}, Journal of Physics: Complexity, 1 (2020), p.~035006,
  \url{https://doi.org/10.1088/2632-072X/ aba8e1}.

\bibitem{DeVille.2021b}
{\sc L.~DeVille}, {\em Consensus on simplicial complexes: Results on stability
  and synchronization}, Chaos, 31 (2021), p.~023137,
  \url{https://doi.org/10.1063/5.0037433}.

\bibitem{DeVille.2015}
{\sc L.~DeVille and E.~Lerman}, {\em Modular dynamical systems on networks},
  Journal of the European Mathematical Society, 17 (2015), pp.~2977--3013,
  \url{https://doi.org/10.4171/JEMS/577}.

\bibitem{Diekman.2013}
{\sc C.~O. Diekman, M.~Golubitsky, and Y.~Wang}, {\em Derived patterns in
  binocular rivalry networks}, Journal of mathematical neuroscience, 3 (2013),
  p.~6, \url{https://doi.org/10.1186/2190-8567-3-6}.

\bibitem{Field.2004}
{\sc M.~J. Field}, {\em Combinatorial dynamics}, Dynamical Systems, 19 (2004),
  pp.~217--243, \url{https://doi.org/10.1080/14689360410001729379}.

\bibitem{Gallo.2022}
{\sc L.~Gallo, R.~Muolo, L.~V. Gambuzza, V.~Latora, M.~Frasca, and
  T.~Carletti}, {\em Synchronization induced by directed higher-order
  interactions}, Communications Physics, 5 (2022),
  \url{https://doi.org/10.1038/s42005-022-01040-9}.

\bibitem{Gandhi.2020}
{\sc P.~Gandhi, M.~Golubitsky, C.~Postlethwaite, I.~Stewart, and Y.~Wang}, {\em
  Bifurcations on fully inhomogeneous networks}, SIAM J. Appl. Dyn. Syst., 19
  (2020), pp.~366--411, \url{https://doi.org/10.1137/18M1230736}.

\bibitem{Golubitsky.2004}
{\sc M.~Golubitsky, M.~Nicol, and I.~Stewart}, {\em Some curious phenomena in
  coupled cell networks}, Journal of Nonlinear Science, 14 (2004),
  pp.~207--236, \url{https://doi.org/10.1007/s00332-003-0593-6}.

\bibitem{Golubitsky.2006}
{\sc M.~Golubitsky and I.~Stewart}, {\em Nonlinear dynamics of networks: The
  groupoid formalism}, Bull. Amer. Math. Soc., 43 (2006), pp.~305--365,
  \url{https://doi.org/10.1090/S0273-0979-06-01108-6}.

\bibitem{Golubitsky.2017b}
{\sc M.~Golubitsky and I.~Stewart}, {\em Homeostasis, singularities, and
  networks}, Journal of Mathematical Biology, 74 (2017), pp.~387--407,
  \url{https://doi.org/10.1007/s00285-016-1024-2}.

\bibitem{Golubitsky.2023}
{\sc M.~Golubitsky and I.~Stewart}, {\em Dynamics and bifurcation in networks:
  Theory and applications of coupled differential equations}, vol.~185 of Other
  titles in applied mathematics, {Society for Industrial and Applied
  Mathematics}, Philadelphia, 2023,
  \url{https://doi.org/10.1137/1.9781611977332}.

\bibitem{Golubitsky.1999c}
{\sc M.~Golubitsky, I.~Stewart, P.-L. Buono, and J.~J. Collins}, {\em Symmetry
  in locomotor central pattern generators and animal gaits}, Nature, 401
  (1999), pp.~693--695, \url{https://doi.org/10.1038/44416}.

\bibitem{Golubitsky.2005}
{\sc M.~Golubitsky, I.~Stewart, and A.~T{\"o}r{\"o}k}, {\em Patterns of
  synchrony in coupled cell networks with multiple arrows}, SIAM J. Appl. Dyn.
  Syst., 4 (2005), pp.~78--100, \url{https://doi.org/10.1137/040612634}.

\bibitem{vonderGracht.2022}
{\sc {\noopsort{Gracht}}{{von der Gracht}, S{\"o}ren}, E.~Nijholt, and
  B.~Rink}, {\em Amplified steady state bifurcations in feedforward networks},
  Nonlinearity, 35 (2022), pp.~2073--2120,
  \url{https://doi.org/10.1088/1361-6544/ac5463}.

\bibitem{Kamei.2009}
{\sc H.~Kamei}, {\em Construction of lattices of balanced equivalence relations
  for regular homogeneous networks using lattice generators and lattice
  indices}, International Journal of Bifurcation and Chaos, 19 (2009),
  pp.~3691--3705, \url{https://doi.org/10.1142/S0218127409025067}.

\bibitem{Levine.2017}
{\sc J.~M. Levine, J.~Bascompte, P.~B. Adler, and S.~Allesina}, {\em Beyond
  pairwise mechanisms of species coexistence in complex communities}, Nature,
  546 (2017), pp.~56--64, \url{https://doi.org/10.1038/nature22898}.

\bibitem{Mulas.2020}
{\sc R.~Mulas, C.~Kuehn, and J.~Jost}, {\em Coupled dynamics on hypergraphs:
  Master stability of steady states and synchronization}, Physical Review E,
  101 (2020), p.~062313, \url{https://doi.org/10.1103/PhysRevE.101.062313}.

\bibitem{Neuhauser.2022}
{\sc L.~Neuh{\"a}user, R.~Lambiotte, and M.~T. Schaub}, {\em Consensus dynamics
  and opinion formation on hypergraphs}, in Higher-Order Systems, F.~Battiston
  and G.~Petri, eds., Understanding Complex Systems, {Springer International
  Publishing}, Cham, 2022, pp.~347--376,
  \url{https://doi.org/10.1007/978-3-030-91374-8_14}.

\bibitem{Nijholt.2022c}
{\sc E.~Nijholt and L.~DeVille}, {\em Dynamical systems defined on simplicial
  complexes: Symmetries, conjugacies, and invariant subspaces}, Chaos
  (Woodbury, N.Y.), 32 (2022), p.~093131,
  \url{https://doi.org/10.1063/5.0093842}.

\bibitem{Nijholt.2022b}
{\sc E.~Nijholt, J.~L. Ocampo-Espindola, D.~Eroglu, I.~Z. Kiss, and
  T.~Pereira}, {\em Emergent hypernetworks in weakly coupled oscillators},
  Nature communications, 13 (2022), p.~4849,
  \url{https://doi.org/10.1038/s41467-022-32282-4}.

\bibitem{Nijholt.2016}
{\sc E.~Nijholt, B.~Rink, and J.~Sanders}, {\em Graph fibrations and symmetries
  of network dynamics}, J. Differential Equations, 261 (2016), pp.~4861--4896,
  \url{https://doi.org/10.1016/j.jde.2016.07.013}.

\bibitem{nijholt2019center}
{\sc E.~Nijholt, B.~Rink, and J.~Sanders}, {\em Center manifolds of coupled
  cell networks}, SIAM Review, 61 (2019), pp.~121--155,
  \url{https://doi.org/10.1137/18M1219977}.

\bibitem{Nijholt.2020}
{\sc E.~Nijholt, B.~Rink, and S.~Schwenker}, {\em Quiver representations and
  dimension reduction in dynamical systems}, SIAM J. Appl. Dyn. Syst., 19
  (2020), pp.~2428--2468, \url{https://doi.org/10.1137/20M1345670}.

\bibitem{Porter.2020}
{\sc M.~A. Porter}, {\em Nonlinearity + networks: A 2020 vision}, in Emerging
  Frontiers in Nonlinear Science, P.~G. Kevrekidis, J.~Cuevas-Maraver, and
  A.~Saxena, eds., vol.~32 of Nonlinear Systems and Complexity, {Springer
  International Publishing}, Cham, 2020, pp.~131--159,
  \url{https://doi.org/10.1007/978-3-030-44992-6_6}.

\bibitem{Salova.2021b}
{\sc A.~Salova and R.~M. D'Souza}, {\em Analyzing states beyond full
  synchronization on hypergraphs requires methods beyond projected networks},
  July 2021, \url{https://arxiv.org/abs/2107.13712}.

\bibitem{Salova.2021}
{\sc A.~Salova and R.~M. D'Souza}, {\em Cluster synchronization on
  hypergraphs}, Mar. 2022, \url{https://arxiv.org/abs/2101.05464}.

\bibitem{Soares.2018}
{\sc P.~Soares}, {\em The lifting bifurcation problem on feed-forward
  networks}, Nonlinearity, 31 (2018), pp.~5500--5535,
  \url{https://doi.org/10.1088/1361-6544/ aae1d0}.

\bibitem{Stewart.2003}
{\sc I.~Stewart, M.~Golubitsky, and M.~Pivato}, {\em Symmetry groupoids and
  patterns of synchrony in coupled cell networks}, SIAM J. Appl. Dyn. Syst., 2
  (2003), pp.~609--646, \url{https://doi.org/10.1137/S1111111103419896}.

\bibitem{Torres.2021}
{\sc L.~Torres, A.~S. Blevins, D.~Bassett, and T.~Eliassi-Rad}, {\em The why,
  how, and when of representations for complex systems}, SIAM Review, 63
  (2021), pp.~435--485, \url{https://doi.org/10.1137/20M1355896}.

\end{thebibliography}

\end{document}